\documentclass[11pt]{extarticle}
\usepackage[english]{babel}
\usepackage{graphicx}
\usepackage{framed}
\usepackage[top=1 in,bottom=1in, left=1 in, right=1 in]{geometry}

\usepackage[utf8]{inputenc} % allow utf-8 input
\usepackage[T1]{fontenc}    % use 8-bit T1 fonts
\usepackage{hyperref}       % hyperlinks
\usepackage{url}            % simple URL typesetting
\usepackage{booktabs}       % professional-quality tables
\usepackage{amsfonts}       % blackboard math symbols
\usepackage{amsmath}
\usepackage{amsthm}
\usepackage{amssymb}
\usepackage{mathtools}
\usepackage{commath}
\usepackage{mathrsfs}
\usepackage{bbm}
\usepackage{nicefrac}       % compact symbols for 1/2, etc.
\usepackage{microtype}      % microtypography
\usepackage{lipsum}
\usepackage{graphicx}
\graphicspath{ {./images/} }
\usepackage{xcolor}
\usepackage{appendix}
\usepackage{enumitem}   
\usepackage{breqn}

\theoremstyle{plain}
\newtheorem{theorem}{Theorem}
\newtheorem{proposition}{Proposition}
\newtheorem{lemma}{Lemma}
\newtheorem{corollary}{Corollary}
\newtheorem{assumption}{Assumption}

\theoremstyle{definition}
\newtheorem*{definition}{Definition}
\newtheorem*{remark}{\textit{Remark}}

\newcommand{\R}{\mathbb{R}}

\newcommand{\C}{\mathcal{C}}

\newcommand{\h}{\mathfrak{h}}
\title{Long time behavior of an age and leaky memory-structured neuronal population equation}
\author{Claudia Fonte\thanks{CEREMADE, Université Paris Dauphine-PSL, 75016 Paris, France (fonte@ceremade.dauphine.fr)} \and Valentin Schmutz\thanks{Brain Mind Institute, École Polytechnique Fédérale de Lausanne, 1015 Lausanne, Switzerland (valentin.schmutz@epfl.ch)}}

\begin{document}
\date{}
\maketitle

\begin{abstract}
   We study the asymptotic stability of a two-dimensional mean-field equation, which takes the form of a nonlocal transport equation and generalizes the time-elapsed neuron network model by the inclusion of a leaky memory variable. This additional variable can represent a slow fatigue mechanism, like spike frequency adaptation or short-term synaptic depression. Even though two-dimensional models are known to have emergent behaviors, like population bursts, which are not observed in standard one-dimensional models, we show that in the weak connectivity regime, two-dimensional models behave like one-dimensional models, i.e. they relax to a unique stationary state. 
   
The proof is based on an application of Harris' ergodic theorem and a perturbation argument, adapted to the case of a multidimensional equation with delays.

\vspace{0.5cm}

\end{abstract}

\tableofcontents

\section{Introduction}
Multidimensional mean-field models in theoretical neuroscience are challenging to analyse \cite{Ric09,VellLin19,BeiOst19,MusGer19} but their study is a necessary step towards understanding how multiple timescales present at the single-neuron level \cite{PozNau13,TeeIye18} affect the dynamics of large networks of neurons.

One-dimensional mean-field equations for populations of spiking neurons with deterministic drift with stochastic jumps have been a subject of mathematical studies since the works of Pakdaman, Perthame and Salort \cite{PakPer09,PakPer13,PakPer14}, providing rigorous foundations to earlier works in theoretical neuroscience \cite{WilCow72,GerHem92, Ger95, Ger00}. These population equations correspond to the mean-field limit of large networks of interacting neurons \cite{DemGal15,FouLoe16,Che17}. However, they are derived from spiking neuron models that are of the `renewal' type (with the exception of \cite{PakPer14}), which means that, while they capture the effect of neuronal refractoriness, they neglect slower neuronal timescales, like those of spike frequency adaptation and short-term synaptic plasticity. 

To take into account slow neuronal timescales, state-of-the-art phenomenological spiking neuron models are multidimensional \cite{KobTsu09,TeeIye18} or kernel-based \cite{TruEde05,PilShl08,PozMen15} (and see \cite[Ch.~6.4]{GerKis14}). In the following, we consider a class of neuron models that characterize neuronal refractoriness by an `age' variable (the time elapsed since last spike) and effects of spike frequency adaptation or short-term synaptic plasticity by a `leaky memory' variables. For this class of neuron models, the mean-field limit is characterized by a multidimensional transport equation with a nonlocal boundary condition \cite{Sch20}. In this work, we study the long time behavior of the solutions to the equation proposed in \cite{Sch20}, in the two-dimensional case.

\subsection{The age- and leaky memory-structured model}
The population model we consider describes the evolution of a density $\rho_t$ over the state-space $(a,m)\in\R_+\times\R_+^*$, where $a$  and $m$ are the `age' and `leaky memory' variables of the neuron, and $\rho_t(a,m)$ represents the density of neurons in state $(a,m)$ at time $t$. 

The nonlinear evolution problem for the density $\rho_t$, for the initial datum $u_0$, writes
\begin{subequations} \label{eq:PDE}
\begin{align} 
    &\partial_t\rho_t+\nabla\cdot (b\rho_t)= - f(a,m,\varepsilon x_t)\rho_t , \label{eq:continuity}\\
    &\rho_t(0,m) = \mathbbm{1}_{m>\gamma(0)}\left|(\gamma^{-1})'(m)\right|\int_0^\infty f(a,\gamma^{-1}(m),\varepsilon x_t)\rho_t(a,\gamma^{-1}(m))da, \label{eq:rho_evol_border} \\
    &x_t = \int_0^t \int_{0}^\infty\int_0^\infty h(t-s,a,m)f(a,m,\varepsilon x_s)\rho_s(a,m)dadmds,\label{eq:rho_evol_x}\\
    &\rho_0 = u_0.
\end{align}
\end{subequations}
The dynamics of the model can be decomposed in three elements: ($i$) the behavior of neurons between spikes, ($ii$) the spike-triggered jumps and ($iii$) the interaction between neurons, which we discuss in turn.

($i$) Between spikes, neurons are transported along the vector field $b(a,m)=(1,-\lambda m)$, with $\lambda > 0$ ($\nabla \cdot$ denotes the divergence operator over the state-space).

($ii$) Neurons spike at a rate $f(a,m,\varepsilon x_t)$, where $f: \R_+\times\R_+^*\times\R \to \R_+$ is the `firing rate function' corresponding to the stochastic intensity of the spike generation process and $\varepsilon \in \R$ is the connection strength. When a neuron spikes, its age $a$ is reset to $0$ and its leaky memory variable $m$ jumps to $\gamma(m)$, where $\gamma:\R_+\to\R_+^*$ is the `jump mapping' and is assumed to be a strictly increasing $\mathcal{C}^1-$diffeomorphism. As a consequence, the border condition \eqref{eq:rho_evol_border} has a simple interpretation: the density of neurons in state $(0,m)$ at time $t$ is equal to the marginal density of those neurons that have their leaky memory variable in state $\gamma^{-1}(m)$ and spike at time $t$. The indicator function $\mathbbm{1}_{m>\gamma(0)}$ reflects the fact that $m$ is always strictly positive and the term $\left|(\gamma^{-1})'(m)\right|$ is necessary to guarantee the conservation of the total mass of neurons. Indeed, formally,
\begin{align*} \partial_t\int \rho_t&=\int \mathbbm{1}_{m>\gamma(0)}\left|(\gamma^{-1})'(m)\right|\int_0^\infty f(a,\gamma^{-1}(m),\varepsilon x_t)\rho_t(a,\gamma^{-1}(m))dadm-\int f(a,m,\varepsilon x_t)\rho_t\\
&=0,
\end{align*} 
by a change of variable.

($iii$) Neurons interact through the `total postsynaptic potential' $x_t$, which integrates the past spiking activity of the population, filtered by the `interaction function' $h:\R_+\times\R_+\times\R_+^*\to\R$. $x_t$, weighted by the connection strength $\varepsilon\in\R$, influences the firing rate $f$. If we write $N(t)$ for the mean firing rate
\begin{equation*}
N(t):=\int_{0}^\infty\int_0^\infty f(a,m,\varepsilon x_t)\rho_t(a,m)dadm,
\end{equation*}
and if we take $h$ independent of $a$ and $m$, then $x_t$ takes the form
\begin{equation*}
     x_t=\int_0^t h(t-s)N(s)ds,
\end{equation*}
where now $h$ is simply a delay kernel, as in \cite{GerHem92,Ger95,Ger00,PakPer09}. In our formulation, $h$ in Eq.~\eqref{eq:rho_evol_x} allows to model more general interactions. For example, in Sec.~\ref{sec:synaptic_dep}, we show that by choosing $h(t,a,m)=\hat{h}(t)(1-m)$, we can include the effects of a classical short-term synaptic plasticity model \cite{TsoPaw98}.

\subsection{Motivation} 

The model~\eqref{eq:PDE} extends the time elapsed neuron network model  \cite{PakPer09} (see also \cite{Ger95,Ger00}) by the addition of a leaky memory variable which can accumulate over spikes (as opposed to the age variable which is reset to $0$ at each spike) and hence introduces a slow timescale in the population dynamics. Such a slow timescale is typically used to account for some form of fatigue mechanism, which can act on the spiking activity (spike frequency adaptation) or on synaptic transmission (short-term synaptic depression). Slow fatigue at the single neuron level can lead to nontrivial emergent behaviors at the population level, like population bursts \cite{VreHan01, GigMat07, GasSch20} (see Fig.~1), which have not been observed in the age- or voltage-structured models of \cite{PakPer09} and \cite{DemGal15} (but see \cite{PakPer14}). 
Even though some population equations have been successfully used in the computational neuroscience literature to study emergent behaviors in networks of neurons with fatigue, these population equations were obtained at the cost of a timescale separation approximation \cite{GigMat07, GasSch20} or a `mixing' assumption \cite{NauGer12, SchDeg17}, making them inexact. In contrast, the model~\eqref{eq:PDE} is the exact mean-field limit \cite{Sch20} for spiking neuron models with spike-frequency adaptation or short-term synaptic depression, as we discuss now. 

\subsubsection{Spike frequency adaptation}
The recent spike history of a neuron can modulate its firing rate $f$, leading to spike frequency adaptation \cite{BenHer03}. If $h$ is independent of $a$ and $m$ and if $\gamma(m) = m + \hat{\Gamma}$, for a fixed $\hat{\Gamma}$>0, \eqref{eq:PDE} becomes
\begin{subequations} \label{eq:PDE_aSRM0}
\begin{align} 
    &\partial_t\rho_t+\nabla\cdot (b\rho_t)= - f(a,m,\varepsilon x_t)\rho_t, \\
    &\rho_t(0,m) =\mathbbm{1}_{m>\hat{\Gamma}}\int_0^\infty f(a,m-\hat{\Gamma},\varepsilon x_t)\rho_t(a,m - \hat{\Gamma})da,  \\
    &x_t = \int_0^t  h(t-s)\int_{0}^\infty\int_0^\infty f(a,m,\varepsilon x_s)\rho_s(a,m)dadmds,\\
    &\rho_0 = u_0.
\end{align}
If $\eta:\R_+\to\R$ is a bounded function such that $\lim_{a\to+\infty}\eta(a)=0$ ($\eta$ is the `refractory kernel' \cite[Sec.~9.3]{GerKis14}), we can define $f$ more explicitly:
\begin{equation} \label{eq:f_aSRM0}
    f(a,m,\varepsilon x_t):=\hat{f}(\eta(a)-m+\varepsilon x_t),
\end{equation}
\end{subequations}
where $\hat{f}:\R\to\R_+$ is typically a non-decreasing function. Since $m$ makes jumps of size $\hat{\Gamma}>0$ at each spike and decays exponentially at rate $\lambda$ between spikes, $m$ accumulates over spikes, which decreases the firing rate $f$ (Eq.~\eqref{eq:f_aSRM0}), leading to spike frequency adaptation \cite{BenHer03}. More specifically, Eq.~\eqref{eq:PDE_aSRM0} is a population equation for adaptive $\text{SRM}_0$ (Spike Response Model) neurons \cite{JolRau06, GerKis14}.

Populations of spiking neurons with spike frequency adaptation exhibit self-sustained population bursts when the connectivity strength is sufficiently strong \cite{VreHan01, GigMat07, GasSch20}. In Fig.~1, we show simulations of \eqref{eq:PDE_aSRM0} for two different connectivity strengths $\varepsilon$. For large $\varepsilon$, we observe self-sustained bursts, whereas for small $\varepsilon$, we observe relaxation to a stationary state.  
For comparison, in the Appendix, we show similar simulations for the time elapsed neuron network model \cite{PakPer09}, where, as expected, we only observe self-sustained oscillations or relaxation to a stationary state.
\begin{figure} \label{fig:1}
  \begin{center}
\includegraphics[scale=0.65]{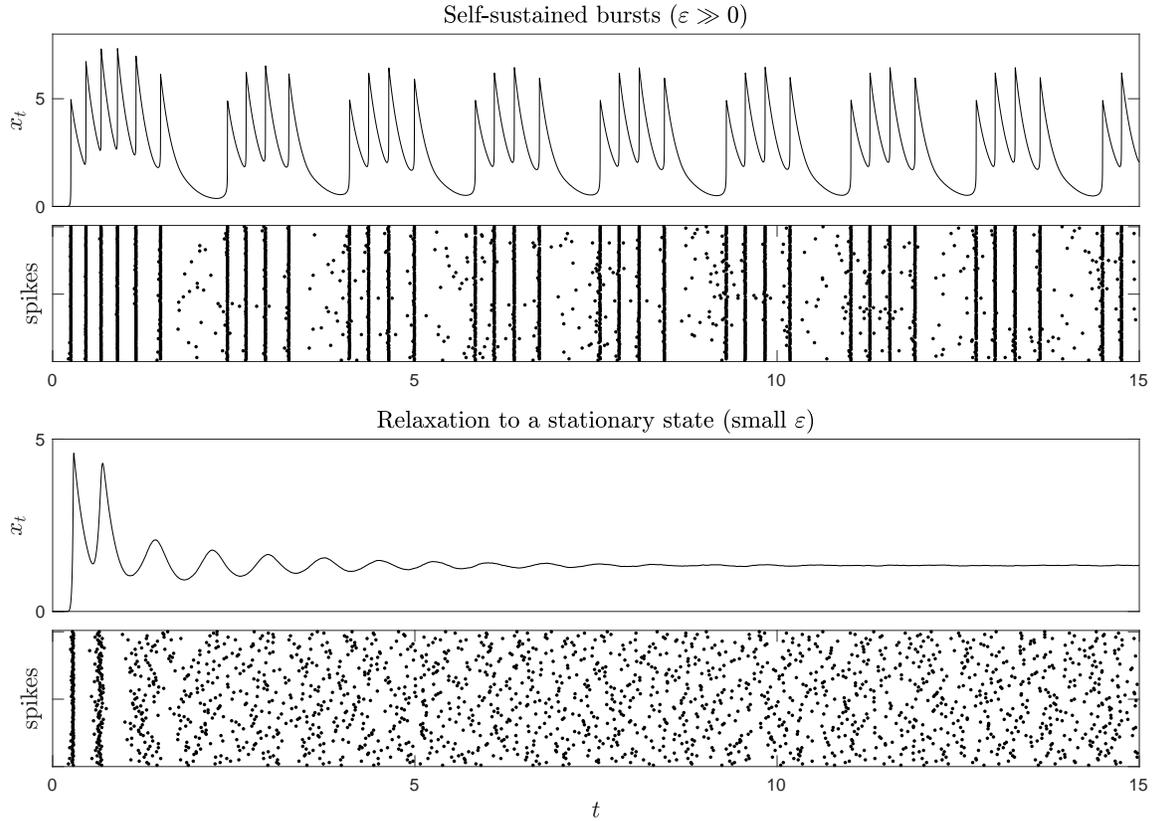}
\caption{\textbf{Depending on the connectivity strength $\varepsilon$, a population of adaptive $\text{SRM}_0$ can exhibit self-sustained bursts ($\varepsilon \gg 0$) or relaxation to a stationary state (small $\varepsilon$)}. We show simulations of a network of $5\cdot 10^5$ adaptive $\text{SRM}_0$ neurons, approximating the mean-field limit Eq.~\eqref{eq:PDE_aSRM0}, with identical parameters (except for $\varepsilon$) and identical initial conditions. The raster plots below the plots for the time-evolution of the total postsynaptic potential $x_t$ represent the spikes of $100$ randomly selected neurons.}
  \end{center}
\end{figure}

\subsubsection{Short-term synaptic depression}\label{sec:synaptic_dep}
The recent spike history of a presynaptic neuron can modulate the synaptic transmission, leading to short-term synaptic plasticity \cite{ZucReg02}. We will consider here the case of depressive synapses and use the model of \cite{TsoPaw98} (with a change of variable for convenience). In this case, the state-space is $(a,m)\in\R_+\times]0,1[$. Taking $f$ independent of $m$, and choosing $h$ and $\gamma$ of the the form $h(t,a,m):=\hat{h}(t)(1-m)$ and $\gamma(m):= 1-\upsilon+\upsilon m$ for a fixed $\upsilon\in]0,1[$, \eqref{eq:PDE} becomes
\begin{subequations} \label{eq:PDE_STD}
\begin{align} 
    &\partial_t\rho_t+\nabla\cdot (b\rho_t)= - f(a,m,\varepsilon x_t)\rho_t,  \\
    &\rho_t(0,m) = \mathbbm{1}_{m>\gamma(0)}\frac{1}{\upsilon}\int_0^\infty f(a,\varepsilon x_t)\rho_t(a,\gamma^{-1}(m))da, \label{eq:PDE_STD_border}\\
    &x_t = \int_0^t \hat{h}(t-s)\int_{0}^1\int_0^\infty (1-m)f(a,\varepsilon x_s)\rho_s(a,m)dadmds, \label{eq:PDE_STD_x}\\
    &\rho_0 = u_0.
\end{align}
\end{subequations}
Note that the term $\frac{1}{\upsilon}$ on the RHS of Eq.~\eqref{eq:PDE_STD_border} simply comes from the fact that $|(\gamma^{-1})'(m)| = \frac{1}{\upsilon}$, for all $m\in]0,1[$.
Here, at each spike, $m$ makes strictly positive jumps which size tends to $0$ as $m$ tends to $1$ (since $\gamma(1)=1$) and decays exponentially at rate $\lambda$ between spikes. If $m$ is close to $1$, synaptic transmission is weak because of the factor $(1-m)$ in Eq.~\eqref{eq:PDE_STD_x}.

As observed in \cite{RomAmi06}, the stationary state of populations of neurons with short-term synaptic plasticity can be described by a simple formula, which we prove in Sec.~\ref{sec:formula_STD}. 

\subsection{Assumptions and main results}
The main result of this work is the exponential stability of~\eqref{eq:PDE} in the weak connectivity regime (\textbf{Theorem~\ref{theorem:stability_nonlinear}}) -- or, more explicitly, there exists $\varepsilon^{**}>0$ such that \eqref{eq:PDE} is exponentially stable for all connectivity strength $\varepsilon \in ]-\varepsilon^{**},+\varepsilon^{**}[$. Before proving the exponential stability, we first establish the well-posedness of~\eqref{eq:PDE} in the appropriate function space (\textbf{Theorem~\ref{theorem:wellposedness}}) and show that stationary solutions exist and are unique for sufficiently weak connectivity (\textbf{Theorem~\ref{theorem:stationary_solutions}}).

Here, we study the weak solutions to~\eqref{eq:PDE} for an initial datum in  $L^1_+:=L^1(\R_+\times\R_+^*,\R_+)$ and write $L^1_+(\R_+^*):=L^1(\R_+^*, \R_+)$.
\begin{definition}[Solutions]
$(\rho,x) \in \C(\R_+, L^1_+) \times \C(\R_+)$ is a solution to~\eqref{eq:PDE}, for the initial datum $u_0 \in L^1_+$, if
\begin{subequations}
\begin{equation} 
    x_t = \int_0^t \int_{0}^\infty\int_0^\infty h(t-s,a,m)f(a,m, \varepsilon x_s)\rho_s(a,m)dadm ds, \qquad \qquad \qquad  \forall t \geq 0, \label{eq:definition_x_weak}
\end{equation}
and if for all $\varphi \in \C^\infty_c(\R_+ \times \R_+ \times \R_+^*)$,
\begin{multline}  \label{eq:PDE_weak}
    0=\int_{0}^\infty\int_0^\infty u_0(a,m)\varphi(0,a, m)dadm \\
    + \int_0^\infty\int_{0}^\infty\int_0^\infty \rho_t(a,m) \Big\{[\partial_t + \partial_a - \lambda m\partial_m] \varphi + (\varphi(t,0, \gamma(m))-\varphi(t,a, m))f(a,m,\varepsilon x_t)\Big\}dadm dt.
\end{multline}
\end{subequations}
\end{definition}

To prove the well-posedness of~\eqref{eq:PDE}, we need some simple assumptions of the firing rate function $f$ and the interaction function $h$:

%To prove the well-posedness of Eq.~\eqref{eq:PDE} and throughout this paper, we assume:

\begin{assumption}\label{assumption:functions} 
$f$ is bounded and $L_f$-Lipschitz, i.e. 
    \begin{equation*}
        |f(a, m, x) - f(a^*, m^*, x^*)| \leq L_f(|a - a^*| + |m - m^*| + |x - x^*|),
    \end{equation*}
and $h$ is bounded and continuous.

\end{assumption}

Since we want to apply Harris' theorem, the well-posedness in $L^1$ (which is treated in \cite{Sch20}) is not enough and we need the well-posedness in a weighted $L^1$ space (where the weight satisfies a Lyapunov condition \cite{MeyTwe93}) with a global-in-time estimate in the weighted $L^1$ norm. 

Using the weight function
\begin{equation*}
    w: \R_+ \times \R_+ \to [1,\infty), \quad (a,m) \mapsto 1+m,
\end{equation*}
we define the function space
\begin{equation*}
    L^1_+(w) := \left\{g\in L^1(\R_+\times\R_+^*,\R_+) \:\Big|\: \norm{g}_{L^1(w)}:=\int_0^\infty \int_0^\infty g(a,m)w(a,m)dadm < \infty\right\}.
\end{equation*}
To obtain a global-in-time estimate in the $L^1_+(w)$ norm, we further need that the jump sizes of $\gamma$ are bounded:

\begin{assumption}\label{assumption:w-estimate}
There exists a bounded function $\Gamma : \R_+^*\to\R_+^*$ such that for all $m\in\R_+^*$, $\gamma(m) = m + \Gamma(m)$. 
\end{assumption}

\begin{theorem}[Well-posedness] \label{theorem:wellposedness}
Grant Assumption~\ref{assumption:functions}. For any initial datum $u_0 \in L^1_+$, there exists a unique weak solution $(\rho,x)$ to~\eqref{eq:PDE}. This solution satisfies
\begin{enumerate}[label=(\roman*)]
    \item \textnormal{($L^1$-stability)} $\quad\norm{\rho_t}_{L^1} = \norm{u_0}_{L^1}$, $\forall t>0$,
    \item \textnormal{(Global bound in $L^1_+(w)$)} if, in addition, Assumption \ref{assumption:w-estimate} holds and $u_0 \in L^1_+(w)$, then
\begin{equation}\label{eq:L_w_estimate}
\forall t>0, \qquad \norm{\rho_t}_{L^1(w)}\leq \norm{u_0}_{L^1(w)}e^{-\alpha t} + \frac{b}{\alpha}(1-e^{-\alpha t}),
\end{equation}
for some constants $\alpha>0$ and $b\in\R$.
\end{enumerate}
\end{theorem}

In contrast to \cite{Sch20}, the well-posedness proof presented here does not involve any probabilistic argument. The proof consists of two consecutive applications of Banach's fixed-point theorem, where a first fixed-point gives the unique solution to a linearized version of~\eqref{eq:PDE} which is then used in a second fixed-point treating the nonlinearity of~\eqref{eq:PDE}.

The second step towards the exponential stability proof is the study of the existence and uniqueness of the stationary solutions to~\eqref{eq:PDE}. For this step, we require:
\begin{assumption}\label{assumption:long-time}
{\color{white}nothing}
\begin{enumerate}[label=(\roman*)]
\item There exists $\Delta_{\text{abs}}>0$ and $\sigma>0$ such that 
\begin{equation*}
f(a,m,x)\geq \sigma, \qquad \forall (a,m,x)\in[\Delta_{\text{abs}},+\infty[\times\R_+^*\times \R.
\end{equation*}
\item There exists $C_\gamma\in]0,1]$ such that $C_\gamma\leq \gamma'\leq 1$.
\item $\bar{h}(a,m)=\int_0^\infty h(t,a,m)dt$ is bounded.
%\item There exists $\Gamma \in C_b^1(\R_+^+,\R_+^*)$ with $\Gamma'\leq 0$ such that $\gamma(m) = m + \Gamma(m), \forall m\in\R_+^*$. 
\end{enumerate}
\end{assumption}

The first point of Assumption~\ref{assumption:long-time} sets a lower bound on the firing rate function $f$ for any $a\geq \Delta_{abs}$ and hence allows for an absolute refractory period $\Delta_{\text{abs}}>0$, i.e. a period of time following a spike during which $f=0$ (which is an important neurodynamical feature \cite[Sec.~1.1]{GerKis14}). This assumption is also used in \cite{CanYol19}. 

In the second point of Assumption~\ref{assumption:long-time}, the lower bound $0<C_\gamma\leq\gamma^\prime$ guarantees that $\gamma$ is strictly increasing, which reflects the idea that $m$ is a `leaky memory' variable of the past neuronal activity. On the other hand, the upper bound $\gamma^\prime\leq 1$, which can be rewritten in terms of the jump size function $\Gamma$ as $\Gamma^{\prime}\leq 0$, prevents the variable $m$ from growing too fast and allows for a potential saturation of the memory, as in the example with short-term synaptic plasticity~\eqref{eq:PDE_STD}. The third point of Assumption~\ref{assumption:long-time} reflects the fact that a single spike has a finite impact on the neuron that receives it. 

We emphasize that the two examples shown above, spike frequency adaptation~\eqref{eq:PDE_aSRM0} and short-term synaptic depression~\eqref{eq:PDE_STD}, satisfy Assumption~\ref{assumption:long-time}.

\begin{theorem}[Stationary solutions]\label{theorem:stationary_solutions}
Grant Assumptions~\ref{assumption:functions} -- \ref{assumption:long-time}.
\begin{enumerate}[label=(\roman*)]
    \item There exists a stationary solution to~\eqref{eq:PDE}.
    \item There exists $\varepsilon^*>0$ such that for all $\varepsilon\in]-\varepsilon^*,+\varepsilon^*[$, the stationary solution to~\eqref{eq:PDE} is unique.
\end{enumerate}
\end{theorem}

Over the course of this work, we obtained the existence of the stationary solution by two different approaches. The first approach is based on the Doeblin-Harris method \cite{HaiMat11} and is similar to that of \cite{CanYol19}. First, we show that when $x_t$ is fixed and time-invariant in \eqref{eq:PDE} (neurons are non-interacting), the system satisfies a Harris condition -- this constitutes a key result of this work --, and we can use Harris' theorem to get the stationary solution. Then, we use the Lipschitz continuity the stationary solutions with respect to the fixed $x$ to prove the existence of a stationary solution for arbitrary connectivity strengths $\varepsilon$. Finally, for $\varepsilon$ small enough, we also get the uniqueness of the stationary solution, by Banach's fixed-point theorem. 

The second approach relies on the fact that the stationary solutions solve an integral equation, for which we can show that a solution exists by Schauder's fixed-point theorem. In the process, we get several estimates on the stationary solutions, namely that they are continuous, bounded, and exponentially decaying in $m$. However, this approach does not give uniqueness.

As mentioned above, the application of Harris' theorem requires us to consider solutions in the weighted space $L^1(w)$. However, in the case where the state-space of the leaky memory variable $m$ is bounded, the situation is simpler: we can use Doeblin's theorem in $L^1$. The following assumption guarantees that $m$ stays in a bounded state-space:

\begin{assumption}\label{assumption:compact}
There exists $G>0$ such that for all $m\in\R_+^*$, $\gamma(m)<G$.
\end{assumption}
Note that this assumption is satisfied in the example with short-term synaptic plasticity~\eqref{eq:PDE_STD}, with $G=1$.

Finally, to study the exponential stability of~\eqref{eq:PDE}, we need an exponential decay on $h$:
\begin{assumption}\label{assumption:nonlinear}
There exists $\h,C_h>0$ such that $h(t,a,m)\leq C_h e^{-\h t}$, $\quad\forall(t,a,m)$.
\end{assumption}

By a perturbation argument similar to that of \cite{MisQui18}, we obtain our main result:

\begin{theorem}[Exponential stability in the weak connectivity regime]\label{theorem:stability_nonlinear}
Grant Assumptions~\ref{assumption:functions} -- \ref{assumption:long-time} and \ref{assumption:nonlinear}. For any $W>0$, there exists $\varepsilon^{**}_W>0$ such that for $\varepsilon \in ]-\varepsilon^{**}_W,+\varepsilon^{**}_W[$, there exists $C\geq 1$ and $c_W>0$ such that for all initial data $u_0\in L^1_+(w)$ with $\norm{u_0}_{L^1}=1$ and $\norm{u_0}_{L^1(w)}\leq W$, the solution $(\rho,x)$ to~\eqref{eq:PDE} satisfies
\begin{equation}\label{eq:final}
    \norm{\rho_t - \rho_\infty}_{L^1(w)} + |x_t-x_\infty| \leq C e^{-c_W t}\left(\norm{u_0 - \rho_\infty}_{L^1(w)} + 1 \right), \qquad \forall t\geq 0,
\end{equation}
where $(\rho_\infty,x_\infty)$ is the unique stationary solution given by Theorem~\ref{theorem:stationary_solutions}~$(ii)$.

If, in addition, we grant Assumption~\ref{assumption:compact}, then there exists $\varepsilon^{**}>0$ such that for all $\varepsilon\in]-\varepsilon^{**},+\varepsilon^{**}[$, there exists $C'\geq 1$ and $c>0$ such that for all initial data $u_0\in L^1_+$ with $\norm{u_0}_{L^1}=1$,
\begin{equation}\label{eq:final_doeblin}
    \norm{\rho_t - \rho_\infty}_{L^1} + |x_t-x_\infty| \leq C' e^{-c t}\left(\norm{u_0 - \rho_\infty}_{L^1} + 1 \right), \qquad \forall t\geq 0.
\end{equation}
\end{theorem}
From the neuronal modeling point of view, this result is not surprising: when the connection strength is weak enough, neurons do not synchronize and the population activity converges to a stationary state. This was already proved for simpler one-dimensional models (see below) and the addition of a leaky memory variable carrying the effect of spike frequency adaptation or short-term synaptic plasticity does not change this behavior. 

%\subsection{Related works}
\subsection{Discussion of the methods}

The asymptotic stability of the age-structured model of \cite{PakPer09} in the weak connectivity regime has been studied using entropy methods (assuming that $f$ is a step-function) \cite{PakPer09,PakPer13}, spectral analysis of semigroups in Banach spaces \cite{MisWen18,MisQui18} or Doeblin's theorem \cite{CanYol19}. For the treatment of the strong connectivity regime, we refer to \cite{PakPer09,MisQui18}. 

On the closely related voltage-structured model \cite{DemGal15}, \cite{CorTan20} also proved asymptotic stability in the weak connectivity regime using Laplace transform techniques. For this model, the nonlinear stability of the stationary solutions has been recently studied in \cite{Cor20} (see also \cite{DroVel21}) and can identify Hopf bifurcations \cite{CorTan20b}.  %However, it is known that even in the absence of a representation formula, the asymptotic stability can be proved and this fact has been illustrated on the model with fatigue of \cite{PakPer14} in the same paper, and in \cite{CanYol19} using Doeblin's theorem.

Doeblin's theorem has also been used in \cite{DumGab20} in the case of the `threshold crossing' neuronal population equation of \cite{OmuKni00}. Note that closely related methods have been used by probabilists to study the ergodicity of single neuron models \cite{HopLoe16,DuaLoe19}.

Our approach combines strategies from \cite{MisWen18} and \cite{CanYol19}, even though \cite{CanYol19} uses Doeblin's instead of Harris' theorem. On the one hand, our proof is based on the application of Harris' theorem for the linear problem, which simplifies the proof of \cite{MisWen18}. On the other hand, we use an argument from in \cite{MisWen18} to deal with delay effects, which are not considered in \cite{CanYol19}. Note that our model is two-dimensional (by the addition of the leaky memory variable), whereas the aforementioned works only considered one-dimensional models.

%Our main tool is a version of Harris' ergodic theorem \cite{Har56,HaiMat11}. 

\subsection{Plan of the paper}
The proof of Theorem~\ref{theorem:wellposedness} (Well-posedness) is presented in \textbf{Section~\ref{sec:well-posedness}}. In \textbf{Section~\ref{sec:exponential}}, we prove the exponential stability of \eqref{eq:PDE} in the non-interacting case $\varepsilon = 0$ using Harris' or Doeblin's theorem. The proof of Theorem~\ref{theorem:stationary_solutions} (Stationary solutions) is presented in \textbf{Section~\ref{sec:stationary}} which is divided in three parts: in the first part, we present a proof which uses the exponential stability of the non-interacting case; in the second part, we present an alternative proof for the existence of stationary solutions which does not involve the Doeblin-Harris method; and in the last part, we present a proof for the formula of \cite{RomAmi06} in the case of short-term synaptic plasticity~\eqref{eq:PDE_STD}. Finally, \textbf{Section~\ref{sec:weak}} is dedicated to the proof of Theorem~\ref{theorem:stability_nonlinear} (Exponential stability in the weak connectivity regime).

\section{Well-posedness}\label{sec:well-posedness}

This section is dedicated to the proof of Theorem~\ref{theorem:wellposedness}, which we decompose is several lemmas. First, we verify the \textit{a priori} $L^1$-stability of the solutions to~\eqref{eq:PDE}, a technical result we use later in the proof. Then, we introduce a linearized version of~\eqref{eq:PDE} and show that it is well-posed by an application of Banach's fixed-point theorem. Another Banach's fixed-point argument is used to treat the nonlinearity of~\eqref{eq:PDE} and concludes the proof of the well-posedness in $L^1$. Finally, we prove the global bound in $L^1_+(w)$ (point $(ii)$ of Theorem~\ref{theorem:wellposedness}), which we will use to apply Harris' theorem in the next sections.

\begin{lemma}[{\textit{A priori} $L^1$-stability}]\label{lemma:L1_stability}
Grant Assumption~\ref{assumption:functions}. If $(\rho,x)$ is a weak solution to~\eqref{eq:PDE} for the initial datum $u_0 \in L^1_+$, then 
\begin{equation*}
    \norm{\rho_t}_{L^1} = \norm{u_0}_{L^1}, \qquad \forall t>0.
\end{equation*}
\end{lemma}
\begin{proof}
By a standard cut-off in time argument, we have that for all $T>0$ and for all $\varphi\in\mathcal{C}_c^\infty(\R_+\times\R_+\times\R_+^*)$,
\begin{multline*}
    \int_0^\infty\int_0^\infty \rho_T(a,m)\varphi(T,a,m)dadm - \int_0^\infty\int_0^\infty u_0(a,m)\varphi(0,a,m)dadm = \\
    \int_0^T\int_0^\infty\int_0^\infty \rho_t(a,m) \Big\{[\partial_t + \partial_a - \lambda m\partial_m] \varphi + (\varphi(t,0, \gamma(m))-\varphi(t,a, m))f(a,m,\varepsilon x_t)\Big\}dadm dt.
\end{multline*}
Let $\chi$ be a function in $\mathcal{C}^\infty_c(\R_+ \times \R_+^*, \R_+)$ such that 
\begin{equation*}
    \chi(a,m) = 1, \quad \text{for all} \quad a^2 + m^2 \leq 1.
\end{equation*}
For all $n \in \mathbb{N}^*$, we write $\widetilde{\varphi}^n \in \mathcal{C}^\infty(\R_+ \times \R_+ \times \R_+^*)$ the classical solution to the transport equation
\begin{subequations}
\begin{align} \label{eq:dual_transport}
    &\partial_t \widetilde{\varphi}^n(t,a,m) + \partial_a \widetilde{\varphi}^n(t,a,m) - \lambda m \partial_m \widetilde{\varphi}^n(t,a,m) = 0,  \\
    &\widetilde{\varphi}^n(0, a, m) = \chi(a/n,m/n). \label{eq:chi}
\end{align}
\end{subequations}

Because of the finite speed of propagation of the transport equation, for all $n$, there exists a function $\varphi^n \in \mathcal{C}^\infty_c(\R_+ \times \R_+ \times \R_+^*)$ such that $\varphi^n(t, a, m) = \widetilde{\varphi}^n(t, a, m)$, for all $(t, a, m) \in [0, T] \times \R_+ \times \R^*$. Hence, for all $n\in\mathbb{N}^*$,
\begin{multline*}
    \int_0^\infty\int_0^\infty \varphi^n(T,a,m)\rho_T(a,m)dadm - \int_0^\infty\int_0^\infty \varphi^n(0,a,m)u_0(a,m)dadm = \\
    \int_0^T\int_0^\infty\int_0^\infty\biggr\{\partial_t \varphi^n+ \partial_a \varphi^n-\lambda m \partial_m \varphi^n+\left(\varphi^n(t,0,\gamma(m)) - \varphi^n(t,a,m)\right)f(a,m,\varepsilon x_t) \biggr\}\rho_t(a,m)dadmdt.
\end{multline*}
As $\varphi_n$ is a solution to Eq.~\eqref{eq:dual_transport} on time $[0, T]$, we get
\begin{multline*}
    \int_0^\infty\int_0^\infty \varphi^n(T,a,m)\rho_T(a,m)dadm - \int_0^\infty\int_0^\infty \varphi^n(0,a,m)u_0(a,m)dadm = \\
    \int_0^T\int_0^\infty\int_0^\infty\biggr\{\left(\varphi^n(t,0,\gamma(m)) - \varphi^n(t,a,m)\right)f(a,m,\varepsilon x_t) \biggr\}\rho_t(a,m)dadmdt.
\end{multline*}
For all $(t, a, m) \in [0, T] \times \R_+ \times \R_+^*$, $\varphi^n(t, a, m) \xrightarrow[n \to \infty]{} 1$, since the initial datum tends to $1$ as $n\to\infty$ (Eq.~\eqref{eq:chi}) and by finite speed of propagation. Thus, by dominated convergence, we get
\begin{equation}
    \int_0^\infty\int_0^\infty \rho_T(a,m)dadm - \int_0^\infty\int_0^\infty u_0(a,m)dadm = 0.
\end{equation}
Since $\rho$ is nonnegative, this concludes the proof.
\end{proof}
%As we say above, the proof of the first part of theorem \ref{theorem:wellposedness} relies on two fixed-point arguments, the first ones, in order to prove for all $x\in\mathcal{C}(\R_+)$ the existence and uniqueness of one solution $\rho^x\in\mathcal{C}(\R_+,L^1_+)$ of the linear equation
Lemma~\ref{lemma:L1_stability} will allow us to prove the well-posedness of~\eqref{eq:PDE} by the means of fixed-point arguments. Let us first introduce a linearized version of Eq.~\eqref{eq:PDE}: for all $x\in\mathcal{C}(\R_+)$, we consider the linear evolution problem
\begin{subequations}\label{eq:linearized}
\begin{align}
    &\partial_t\rho_t+ \nabla \cdot (b\rho_t) = - f(a,m,\varepsilon x_t)\rho_t,  \label{eq:lin_continuity}\\
    &\rho_t(0,m) = \mathbbm{1}_{m>\gamma(0)}\left|(\gamma^{-1})'(m)\right|\int_0^\infty f(a,\gamma^{-1}(m),\varepsilon x_t)\rho_t(a,\gamma^{-1}(m))da,  \label{eq:lin_border}\\
    &\rho_0 = u_0.
\end{align}
\end{subequations}
%The second, to obtain a solution of the nonlinear equation from the linearized one.
%The linear equation above will have a special role in the following sections. We study it separately in the next proposition. 

We can see Eq.~\eqref{eq:linearized} as the Kolmogorov forward equation of a time-dependent Makrov process. Indeed, we can rewrite Eqs.~\eqref{eq:lin_continuity} and \eqref{eq:lin_border} as
\begin{equation}
    \partial \rho_t = \mathcal{L}_{t}\rho_t
\end{equation}
where, for all suitable test function $\phi:\R_+\times\R_+^* \to \R$,
\begin{equation}\label{eq:generator}
    \mathcal{L}^*_{t}\phi(a,m) = b(a,m)\cdot \nabla\phi(a,m) + [\phi(0,\gamma(m)) - \phi(a,m)]f(a,m,\varepsilon x_t).
\end{equation}
$\mathcal{L}^*_{t}$ is the time-dependent generator of a piecewise deterministic Markov process with degenerate jumps.

The linearized equation~\eqref{eq:linearized} will play a special role in the following sections and it therefore deserves its own proposition:

%The next proposition concerns its well-posedness.

\begin{proposition}[Well-posedness of the linearized equation~\eqref{eq:linearized}]\label{proposition:well-posedness-linear}
Grant Assumption~\ref{assumption:functions}. For any initial datum $u_0\in L^1_+$ and any $x\in\mathcal{C}(\R_+)$, there exists a unique weak solution $\rho^x \in \C(\R_+,L^1_+)$ to Eq.~\eqref{eq:linearized}. Furthermore, $\rho^x$ satisfies
\begin{enumerate}[label=(\roman*)]
    \item For all $t>0$ and for all $m\in\R_+^*$,
\begin{subequations}
\begin{align*}
    \rho^x_{t}(0,m) &= \mathbbm{1}_{m>\gamma(0)}\left|(\gamma^{-1})'(m)\right|\int_0^\infty f(a,\gamma^{-1}(m),\varepsilon x)\rho^x_{t}(a,\gamma^{-1}(m))da,\\
    \rho^x_t(a,m) &= \begin{cases} 
    u_0(a-t,e^{\lambda t}m)\exp\left(\lambda t-\int_0^t f(a-t+s,e^{\lambda(t-s)}m,\varepsilon x)ds\right) &\text{if}\quad a \geq t, \\
    \rho^x_{t-a}(0,e^{\lambda a}m)\exp\left(\lambda a - \int_{t-a}^t f(a-t+s,e^{\lambda(t-s)}m, \varepsilon x)ds\right) &\text{if}\quad 0 < a < t.
    \end{cases}
\end{align*}
\end{subequations}
    \item For all $t>0$ and for all $\phi\in\C^\infty_c(\R_+\times\R_+^*)$,
\begin{equation}
    \langle \rho^x_t,\phi\rangle = \langle u_0,\phi\rangle + \int_0^t \langle\rho^x_t,\mathcal{L}^*_{x}\phi\rangle ds.
\end{equation}
\end{enumerate}
\end{proposition}

\begin{proof}
Fix $x\in\mathcal{C}(\R_+)$. For all $p\in\C(\R_+,L^1_+(\R_+^*))$ and $u_0\in L^1_+$, we know, from the standard theory of transport equations, that there is a unique weak solution to 
\begin{subequations}
\begin{align*}
    &\partial_t\rho_t+ \nabla\cdot(b\rho_t) = - f(a,m,\varepsilon x_t)\rho_t,  \\
    &\rho_t(0,m) = p_t(m), \\
    &\rho_0 = u_0,
\end{align*}
\end{subequations}
which we denote $\rho^{x,p}$ and is given by the representation formula,
\begin{equation*}
    \rho^{x,p}_t(a,m) := \begin{cases} 
    u_0(a-t,e^{\lambda t}m)\exp\left(\lambda t-\int_0^t f(a-t+s,e^{\lambda(t-s)}m, \varepsilon x_s)ds\right) &\text{if}\quad a \geq t, \\
    p_{t-a}(e^{\lambda a}m)\exp\left(\lambda a - \int_{t-a}^t f(a-t+s,e^{\lambda(t-s)}m, \varepsilon x_s)ds\right) &\text{if}\quad 0 < a < t.
    \end{cases}
\end{equation*}

$\rho^{x,p}$ is in $\C(\R_+,L^1)$, since
\begin{equation*}
    \forall t\in\R_+,\qquad \norm{\rho_t^{x,p}}_{L^1} \leq \norm{u_0}_{L^1} + \int_0^t\norm{p_s}_{L^1}ds.
\end{equation*}

We have 
\begin{equation*}
    \left( \mathbbm{1}_{m>\gamma(0)}\left|(\gamma^{-1})'(m)\right|\int_0^\infty f(a,\gamma^{-1}(m),x)\rho^{x,p}_t(a,\gamma^{-1}(m))da\right)_{(t,m)\in\R_+\times\R_+^*} \in \C(\R_+,L^1_+(\R_+^*))
\end{equation*}
since
\begin{multline*}
    \forall t\in\R_+, \qquad \int_{\gamma(0)}^\infty \left|(\gamma^{-1})'(m)\right|\int_0^\infty f(a,\gamma^{-1}(m),x)\rho^{x,p}_t(a,\gamma^{-1}(m))dadm \\
    \leq \norm{f}_\infty\norm{\rho_t^{x,p}}_{L^1}\leq \norm{f}_\infty\left(\norm{u_0}_{L^1} + \int_0^t\norm{p_s}_{L^1}ds\right).
\end{multline*}
Hence, we can define, for any $T>0$, the operator $\Phi_T^x$:
\begin{align*}
    \C([0,T],L^1_+(\R_+^*)) &\to \C([0,T],L^1_+(\R_+^*)) \\
    p &\mapsto \left( \mathbbm{1}_{m>\gamma(0)}\left|(\gamma^{-1})'(m)\right|\int_0^\infty f(a,\gamma^{-1}(m),x)\rho^{x,p}_t(a,\gamma^{-1}(m))da\right)_{(t,m)\in[0,T]\times\R_+^*}.
\end{align*}

For any $p,q \in \C([0,T], L^1_+(\R_+^*))$, 
\begin{align*}
    \norm{\Phi_T^x(p) - \Phi_T^x(q)}_{\C([0,T],L^1)} &\leq \norm{f}_\infty\sup_{t\in[0,T]}\norm{\rho_t^{x,p}-\rho_t^{x,q}}_{L^1} \\
    &\leq \norm{f}_\infty\int_0^T\norm{p_s-q_s}_{L^1}ds \\
    &\leq T\norm{f}_\infty\norm{p-q}_{\C([0,T], L^1)}. 
\end{align*}
Therefore, if $0<T<\norm{f}_{\infty}^{-1}$, $\Phi^x_T$ is a contraction. By Banach's fixed-point theorem,
there exists a unique $\rho^x\in\C([0,T],L^1_+)$ solving Eq.~\eqref{eq:linearized}. Since the choice of the contracting $T$ does not depend on the initial datum, we can iterate the above argument on successive time intervals of length $T$ and conclude that there exists a unique $\rho^x\in\C(\R_+,L^1_+)$ solving Eq.~\eqref{eq:linearized} for which the formula $(i)$ is satisfied. Then, $(ii)$ follows from a standard cut-off in time argument. 
\end{proof}

Now, we can prove the existence and uniqueness of a solution to the nonlinear problem~\eqref{eq:PDE} by the means of a second application of Banach's fixed-point theorem.

\begin{proof}[Proof of the well-posedness of~\eqref{eq:PDE} in $L^1$]
For any $x\in\mathcal{C}(\R_+)$, we take the $\rho^x$ given by Proposition~\ref{proposition:well-posedness-linear}. We have
\begin{equation*}
    \left(\int_0^t \int_{\R_+\times\R_+^*} h(t-s)f(\varepsilon x_s)\rho_s^x\,dadmds\right)_{t\in\R_+} \in \C(\R_+)
\end{equation*}
since 
\begin{equation*}
    \forall t\in\R_+, \qquad \left|\int_0^t \int_{\R_+\times\R_+^*} h(t-s)f(\varepsilon x_s)\rho_s^x\,dadmds\right| \leq \norm{h}_{\infty}\norm{f}_{\infty}\int_0^t\norm{\rho_s^x}_{L^1}ds.
\end{equation*}
Hence, for any $T>0$, we can define the operator
\begin{align*}
    \Psi_T : \C([0,T]) &\to \C([0,T])\\
    x &\mapsto \left(\int_0^t \int_{\R_+\times\R_+^*} h(t-s)f(\varepsilon x_s)\rho_s^x\,dadmds\right)_{t\in[0,T]}.
\end{align*}
For any $x,y\in\C([0,T])$, we have
\begin{align*}
    \norm{\Psi_T(x)-\Psi_T(y)}_{\C([0,T])} &\leq T\norm{h}_{\infty}\sup_{t\in[0,T]}\int_{\R_+\times\R_+^*} |f(\varepsilon x_t)\rho^x_t - f(\varepsilon y_t)\rho^y_t|\,dadm \\
    &\leq T\norm{h}_\infty\sup_{t\in[0,T]}\left(\varepsilon L_f|x_t-y_t|\;\norm{\rho^x_t}_{L^1} + \norm{f}_\infty\norm{\rho^x_t - \rho^y_t}_{L^1}\right).
\end{align*}
By Grönwall's lemma, $\norm{\rho^x_t}_{L^1} \leq \norm{u_0}_{L^1}\exp(\norm{f}_\infty t)$, since
\begin{equation*}
    \forall t\in[0,T], \qquad \norm{\rho^x_t}_{L^1} \leq \norm{u_0}_{L^1} + \norm{f}_\infty \int_0^t\norm{\rho^x_s}_{L^1}ds. 
\end{equation*}
On the other hand, we have, for all $t\in[0,T]$,
\begin{align*}
    \norm{\rho^x_t-\rho^y_t}_{L^1} &\leq \int_0^t\int_0^\infty \Bigg|\rho^x_s(0,m)\exp\left(-\int_s^t f(u-s,e^{-\lambda(u-s)}\gamma(m),\varepsilon x_u)du\right)\\
    &\qquad\qquad\qquad\qquad\qquad -\rho^y_s(0,m)\exp\left(-\int_s^t f(u-s,e^{-\lambda(u-s)}\gamma(m),\varepsilon y_u)du\right)\Bigg|dmds \\
    &\leq \norm{f}_\infty\int_0^t\norm{\rho^x_s-\rho^y_s}_{L^1}ds + t\varepsilon \norm{f}_\infty  L_f \norm{x-y}_{\mathcal{C}([0,T])}\int_0^t\norm{\rho^x_s}_{L^1}ds.
\end{align*}
Hence, by Grönwall's lemma, for all $t\in[0,T]$,
\begin{align*}
    \norm{\rho^x_t-\rho^y_t}_{L^1} &\leq \varepsilon L_f \norm{u_0}_{L^1}\frac{\left(\exp(\norm{f}_\infty t) - 1\right)^2}{\norm{f}_\infty}\norm{x-y}_{\mathcal{C}([0,T])}. 
\end{align*}
Gathering the bounds, we get
\begin{equation*}
    \norm{\Psi_T(x)-\Psi_T(y)}_{\C([0,T])} \leq T\varepsilon\norm{h}_\infty L_f \norm{u_0}_{L^1}\exp(\norm{f}_\infty T)\left[1 +   \exp(\norm{f}_\infty T)\right]\norm{x-y}_{\mathcal{C}([0,T])}.
\end{equation*}
For $T$ small enough, $\Psi_T$ is a contraction and, by Banach's fixed-point theorem, has a unique fixed-point. Thus, there exists a unique solution $(\rho,x)\in\mathcal{C}([0,T],L^1_+)$. Since, by Lemma~\ref{lemma:L1_stability}, $\norm{\rho_T}_{L^1} = \norm{u_0}_{L^1}$, we can iterate this argument on successive time intervals of length $T$ and conclude that there exists a unique solution in $\mathcal{C}(\R_+,L^1_+)$.
\end{proof}

To conclude the proof of Theorem~\ref{theorem:wellposedness}, it remains to show the estimate Eq.~\eqref{eq:L_w_estimate}. Under Assumption~\ref{assumption:w-estimate}, the weight function 
\begin{equation*}
    w: \R_+ \times \R_+ \to [1,\infty), \qquad (a,m) \mapsto 1+m
\end{equation*}
satisfies $w(a,m)\to\infty$ when $m \to \infty$ and the Lyapunov condition on $m$:
\begin{equation}\label{eq:lyapunov}
    \exists \alpha > 0, b\geq 0 \quad \text{ such that } \quad \mathcal{L}^*_t w \leq -\alpha w + b.  
\end{equation}
Indeed, for all $(t,a,m)\in\R_+\times\R_+\times\R_+^*$,
\begin{align*}
    \mathcal{L}^*_t w(a,m) &= - \lambda m + \Gamma(m)f(a,m,\varepsilon x_t) \leq - \lambda w(a,m) + \lambda +\norm{\Gamma}_\infty\norm{f}_{\infty}.
\end{align*}
Importantly, the constants $\alpha$ and $b$ do not depend on $x$.

\begin{lemma}[Global bound in $L^1_+(w)$] \label{lemma:bounds_l1w}
Grant Assumptions~\ref{assumption:functions} and \ref{assumption:w-estimate}. If the initial datum $u_0$ is in $L^1_+(w)$, then $\rho_t\in L^1_+(w)$ for all $t\geq 0$. Moreover,
\begin{equation}\label{eq:lemma_bound}
\forall t>0, \qquad \norm{\rho_t}_{L^1(w)}\leq \norm{u_0}_{L^1(w)}e^{-\alpha t} + \frac{b}{\alpha}(1-e^{-\alpha t}),
\end{equation}
where the constants $\alpha$ and $b$ are taken from the Lyapunov condition~\eqref{eq:lyapunov}.
\end{lemma}
\begin{proof}
We divide the proof in two steps: first,  we prove that the solution is stable in $L^1_+(w)$ with a weaker and time dependent bound; then, we use this first bound to apply the dominated convergence theorem and obtain Eq.~\eqref{eq:lemma_bound} by Grönwall's lemma.

\textbf{Step 1.} Fix any $T>0$. Let $\chi\in\C^\infty_c(\R_+,\R_+)$ be a non-increasing function such that $\chi(x) = 1$ if $0\leq x\leq 1$ and %$ \chi^\prime \leq C\chi$. 
$\chi(x)=0$ if $x>2$. For all $n\in\mathbb{N}^*$, let us write $\varphi_k(a)\chi_n(m):=\chi(a/k)\chi(m/n)$. For all $n,k$, $w\chi_n\varphi_k \in \C^\infty_c(\R_+\times\R_+,\R_+)$.
Hence, 
\begin{equation*}
    \forall n\in\mathbb{N}^*, \qquad \langle \rho_{T} ,w\chi_n\varphi_k\rangle = \langle u_0,w\chi_n\varphi_k\rangle + \int_0^T \langle\rho_t,\mathcal{L}^*_{x}(w\chi_n\varphi_k)\rangle dt,
\end{equation*}
where
\begin{eqnarray*}
    \mathcal{L}^*_{x}(w\chi_n\varphi_k) &=& \partial_a(w\chi_n\varphi_k) - \lambda m \partial_m(w\chi_n\varphi_k) + \left(w(\gamma(m))\chi_n(\gamma(m))\varphi_k(0) - w\chi_n\varphi_k\right)f\\
    &=& w\chi_n \frac{1}{k} \chi^\prime(a/k)+w\varphi_k\frac{1}{n}\chi^\prime(m/n)- \lambda m w\chi_n + (w(\gamma(m))\chi_n(\gamma(m))\varphi_k - w\chi_n\varphi_k)f.
\end{eqnarray*}
From the $L^1$-stability and the fact that both $w\partial_m \chi_n$ and $w\chi_n$ are bounded, we take the limit $k\to\infty$ with the dominated convergence theorem to obtain
\begin{equation}\label{eq:weakchiformula}
    \langle \rho_{T},w\chi_n\rangle = \langle u_0,w\chi_n\rangle + \int_0^T \left\langle\rho_t,w\frac{1}{n} \chi^{\prime}(m/n)- \lambda m w\chi_n + (w(\gamma(m))\chi_n(\gamma(m))- w\chi_n)f\right\rangle dt.
\end{equation}
From the properties of $\chi$ and $\gamma$, we get
\begin{align*}
    w(0,\gamma(m))\chi_n(\gamma(m)) &\leq w(0,m + \norm{\Gamma}_\infty)\chi_n(\gamma(m))\leq(1+\norm{\Gamma}_\infty)w(a,m)\chi_n(m), \\ 
    \frac{1}{n} w \chi^\prime(m/n) &\leq \dfrac{1+2n}{n}\|\chi^\prime\|_\infty.
\end{align*}
Then, since $f$ is bounded, there exists a constant $C$,  which does not depend on $n$, such that
\begin{equation*}
    \langle \rho_{T},w\chi_n\rangle \leq \langle u_0,w\chi_n\rangle + \int_0^T \langle\rho_t,C(w\chi_n(m)+1)\rangle dt.
\end{equation*}
We can now apply Grönwall's lemma to obtain 
\begin{equation*}
    \langle \rho_{T},w\chi_n\rangle \leq \max( \langle u_0,w\chi_n\rangle C e^{Ct},C).
\end{equation*}
It follows from Fatou's lemma that $\rho_t=S_{T}^{x}u_0\in L_+^1(w)$.

\textbf{Step 2.} To improve the previous estimate, we come back to \eqref{eq:weakchiformula} and use dominated convergence in $n$ (domination being guaranteed by the bound above) to get
\begin{equation*}
    \langle \rho_{T},w\rangle = \langle u_0,w\rangle + \int_0^T \langle\rho_t,\mathcal{L}^*_{x}w\rangle dt.
\end{equation*}
By the Lyapunov condition~\eqref{eq:lyapunov},
\begin{equation*}
    \norm{\rho_{T}}_{L^1(w)} \leq \norm{u_0}_{L^1(w)} - \alpha\int_0^T \norm{\rho_t}_{L^1(w)} + Tb.
\end{equation*}
Finally, Eq.~\eqref{eq:lemma_bound} is obtained by Grönwall's lemma.
\end{proof}

\begin{remark}
Following the same steps as in the proof above, we can show that the bound Eq.~\eqref{eq:lemma_bound} also holds for the linearized equation~\eqref{eq:linearized} and does not depend on $x$ nor the constants $\alpha$ and b.
\end{remark}

\section{Exponential stability in the non-interacting case}\label{sec:exponential}
%As we have mentioned, the linear equation Eq.~\eqref{eq:linearized} plays an important role in the argument of this work. In this section, we prove the long-time convergence of every solution of this equation to the only stationary solution of this dynamics. Later in section \ref{sec:weak}, a perturbation argument is used to demonstrate the stability of the stationary solution of the nonlinear system under the regime of weak connections.

If $x\in\mathcal{C}(\R_+)$ in the linearized equation~\eqref{eq:linearized} is time-invariant, i.e. $x \equiv \tilde{x}$ for some $\tilde{x}\in\R$, then Eq.~\eqref{eq:linearized} can be seen as the dynamics of a non-interacting population of neurons. In this section, we prove the exponential stability in the non-interacting case using Harris' or Doeblin's theorem. This is the key result of this work and will allow us to prove the existence and uniqueness of the stationary solution to \eqref{eq:PDE} (Section~\ref{sec:stationary}) and the exponential convergence to it (Section~\ref{sec:weak}).

For $\tilde{x}\in \R$, $u_0\in L^1$, we denote $\rho^{\tilde{x}}_t$ the unique solution to Eq.~\eqref{eq:linearized} for the initial datum $u_0$ and $x\equiv \tilde{x}$, given by Proposition~\ref{proposition:well-posedness-linear}. We write, using the semigroup notation,
\begin{equation}\label{eq:semigroup}
     S_t^{\tilde{x}}u_0 := \rho_t^{\tilde{x}}, \qquad\forall t\geq0. 
\end{equation}

To show that the Eq.~\eqref{eq:semigroup} is exponentially stable we will use Harris' theorem in the general case or Doeblin's theorem if Assumption~\ref{assumption:compact} is granted. In both cases, the main technical difficulty is to verify the Doeblin minoration condition (Lemma~\ref{lemma:doeblin}) as the jumps of the process described by Eq.~\eqref{eq:generator} are degenerate and the model is two-dimensional.

\begin{lemma}[Doeblin minoration condition]\label{lemma:doeblin}
Grant Assumptions~\ref{assumption:functions} and \ref{assumption:long-time}. Fix any $x\in\R$. For all $R>0$, there exists $T>0$ and a positive non-zero measure $\nu$ such that
\begin{equation}\label{eq:doeblin}
    \forall u_0\in L^1_+, \qquad S^{\tilde{x}}_{T}u_0\geq \nu\int_{\R_+\times]0,R]} u_0\,dadm.
\end{equation}
\end{lemma} 

\begin{proof}
We proceed in two steps. First (\textbf{Step~1}), we choose a time $T>0$ and a rectangle $[0,\bar{a}]\times[\underline{m},\overline{m}]\subset\R_+\times\R_+^*$ (with nonzero Lebesgue measure) and show that the density $S^{\tilde{x}}_{T}u_0\in L^1$ has a lower bound on $[0,\bar{a}]\times[\underline{m},\overline{m}]$ which depends on a Lebesgue integral in $\R_+^2$ involving $u_0$. Then (\textbf{Step~2}), we perform a change of variable to express this lower bound in terms of $\int_{\R_+\times]0,R]} u_0\,dadm$. The proof only relies on the expression of $S^{\tilde{x}}_{t}u_0$ given by the method of characteristics (see Proposition~\ref{proposition:well-posedness-linear}) and this allows treating a typically probabilistic question -- the Doeblin minoration condition -- from a transport point of view. This is possible because $S^{\tilde{x}}_{t}$ is the stochastic (mass-conservative) semigroup of a piecewise deterministic Markov process. 

The constants $\Delta_{\text{abs}}$, $\sigma$ and $C_\gamma$ are taken from Assumption~\ref{assumption:long-time}.

\textbf{Step~1:}

Fix $R>0$. Since $\gamma(e^{-\lambda \Delta_{\text{abs}}}\gamma(0))>\gamma(0)$ and $\gamma(e^{-\lambda t}\gamma(e^{-\lambda \Delta_{\text{abs}}}R)) \to \gamma(0)$ as $t\to\infty$, there exists $\bar{a}>0$ and $T>\bar{a}+\Delta_{\text{abs}}$ such that
\begin{equation}\label{eq:overline_m}
    \underline{m}=:\gamma(e^{-\lambda(T-\bar{a}-\Delta_{\text{abs}})}\gamma(e^{-\lambda \Delta_{\text{abs}}}R))<e^{-\lambda\bar{a}}\gamma(e^{-\lambda \Delta_{\text{abs}}}\gamma(0))=:\overline{m}.
\end{equation}
Eq.~\eqref{eq:overline_m} has the following heuristic interpretation: if we see $S^{\tilde{x}}_{t}$ as the stochastic semigroup of the piecewise deterministic Markov process defined by the generator Eq.~\eqref{eq:generator},  for any initial point $(a_0,m_0)\in \R_+\times]0,R]$ and any landing point $(a,m)\in [0,\bar{a}]\times[\underline{m},\overline{m}]$ at time $T$, there is a `possible' trajectory going from $(a_0,m_0)$ to $(a,m)$, with exactly two jumps (spikes). Since the trajectories of the process are determined by the jump times, we will exploit the fact that these `possible' trajectories correspond to jump times with strictly positive probability density. Below, we take a transport point of view on this probabilistic argument.

For all $(a,m)\in [0,\bar{a}]\times[\underline{m},\overline{m}]$, 
\begin{align*}
    (S^{\tilde{x}}_{T}u_0)(a,m)&\geq\mathbbm{1}_{\{a<T\}}(S^{\tilde{x}}_{T-a}u_0)(0,e^{\lambda a}m)\exp\left(\lambda a - \int_{T-a}^T f(a-T+s, e^{\lambda(T-s)}m,\tilde{x})ds\right) \\
    &\geq\mathbbm{1}_{\{a<T\}}e^{-\norm{f}_\infty T}e^{\lambda a}(S^{\tilde{x}}_{T-a}u_0)(0,e^{\lambda a}m) \\
    &\geq \mathbbm{1}_{\{a<T\}}e^{-\norm{f}_\infty T}\sigma e^{\lambda a}\left|(\gamma^{-1})'(e^{\lambda a}m)\right|\int_{\Delta_{\text{abs}}}^\infty (S^{\tilde{x}}_{T-a}u_0)(a',\gamma^{-1}(e^{\lambda a}m))da'\\
    &= \mathbbm{1}_{\{a<T\}}e^{-\norm{f}_\infty T}\sigma\left|\frac{d}{dm}\gamma^{-1}(e^{\lambda a}m)\right|\int_{\Delta_{\text{abs}}}^\infty \underbrace{(S^{\tilde{x}}_{T-a}u_0)(a',\gamma^{-1}(e^{\lambda a}m))}_{(\star)}da'.\\
\end{align*}
Above, we went back in time to the last jump time $T-a$. Let us notice that $\gamma^{-1}(e^{\lambda a} m)\geq \gamma^{-1}(e^{\lambda a} \underline{m})>0$. We can therefore define
\begin{equation*}
    a^*_{a,m} :=\frac{1}{\lambda}\left(\log \gamma(0) - \log \gamma^{-1}(e^{\lambda a}m)\right).
\end{equation*}
Note that $a^*_{a,m}$ satisfies $\gamma^{-1}(e^{\lambda a^*_{a,m}}\gamma^{-1}(e^{\lambda a}m))=0$. In other words, $a^*_{a,m}$ is the minimal time between the last and second last jumps for a trajectory landing at $(a,m)$ at time $T$. We can easily verify that, by our choice of $\{T, \bar{a}, \underline{m}, \overline{m}\}$, $\Delta_{\text{abs}}\leq a^*_{a,m} < T-a-\Delta_{\text{abs}}$. This guarantees that it is possible to make two jumps in $[0,T]$ and land at $(a,m)$ at time $T$ while respecting the absolute refractoriness of the neuron (i.e. there needs to be a time interval $\geq \Delta_{\text{abs}}$ between jumps). This allows us to go further back in time to the second last jump:

%$m\leq \overline{m}$, we can easily verify that $a^*_{a,m}\geq \Delta_{\text{abs}}$. 

%On the other hand, 
%\begin{multline}\label{eq:doeblin_useful_ineq}
%    \gamma^{-1}\left(e^{\lambda(T-a-\Delta_{\text{abs}})}\gamma^{-1}(e^{\lambda a}m)\right) > \gamma^{-1}\left(e^{\lambda(T-\overline{a}-\Delta_{\text{abs}})}\gamma^{-1}(m)\right) \geq \gamma^{-1}\left(e^{\lambda(T-\overline{a}-\Delta_{\text{abs}})}\gamma^{-1}(\underline{m})\right) \\
%    = e^{-\lambda \Delta_{\text{abs}}}R>0,
%\end{multline}
%whence
%\begin{equation*}
%    \gamma^{-1}\left(e^{\lambda(T-a-\Delta_{\text{abs}})}\gamma^{-1}(e^{\lambda a}m)\right) > \gamma^{-1}\left(e^{\lambda a^*_{a,m}}\gamma^{-1}(e^{\lambda a}m)\right),
%\end{equation*}
%which implies $T-a-\Delta_{\text{abs}}>a^*_{a,m}$. In summary, we have
%\begin{equation}
%    \Delta_{\text{abs}}\leq a^*_{a,m} < T-a-\Delta_{\text{abs}}.
%\end{equation}
For all $a'\in[a^*_{a,m},T-a-\Delta_{\text{abs}}]$,
\begin{multline*}
    (\star) \geq \mathbbm{1}_{\{a'<T-a\}}e^{-\norm{f}_\infty T}\sigma\left|(\gamma^{-1})'(e^{\lambda a'}\gamma^{-1}(e^{\lambda a}m))\right|e^{\lambda a'}\\
    \int_{\Delta_{\text{abs}}}^\infty\underbrace{(S^{\tilde{x}}_{T-a-a'}u_0)(a'', \gamma^{-1}(e^{\lambda a'}\gamma^{-1}(e^{\lambda a}m)))}_{(\star\star)}da''.
\end{multline*}
Then, we can go further back to time $0$ to get $u_0$:
\begin{multline*}
    (\star\star) \geq \mathbbm{1}_{\{a''\geq T-a-a'\}}e^{-\norm{f}_\infty T}e^{\lambda(T-a-a')}u_0(a''-(T-a-a'),e^{\lambda(T-a-a')}\gamma^{-1}(e^{\lambda a'}\gamma^{-1}(e^{\lambda a}m))).
\end{multline*}
Putting all the lower bounds together, we get
\begin{multline*}
    (S^{\tilde{x}}_{T}u_0)(a,m)\geq\mathbbm{1}_{\{a<T\}}e^{-3\norm{f}_\infty T}\sigma^2\\
    \int_{a^*_{a,m}}^{T-a-\Delta_{\text{abs}}}\int_{T-a-a'}^\infty \left|\frac{d}{dm}e^{\lambda(T-a-a')}\gamma^{-1}(e^{\lambda a'}\gamma^{-1}(e^{\lambda a}m))\right|\\
    u_0(a'' - (T-a-a'), e^{\lambda(T-a-a')}\gamma^{-1}(e^{\lambda a'}\gamma^{-1}(e^{\lambda a}m)))da''da'.
\end{multline*}
Since $\gamma' \leq 1$ (Assumption~\ref{assumption:long-time}), 
\begin{equation*}
    \left|\frac{d}{dm}e^{\lambda(T-a-a')}\gamma^{-1}(e^{\lambda a'}\gamma^{-1}(e^{\lambda a}m))\right|\geq e^{\lambda T}.
\end{equation*}
Thus,
\begin{equation}\label{eq:lower_bound_rho_T}
    (S^{\tilde{x}}_{T}u_0)(a,m)\geq \mathbbm{1}_{\{a<T\}}e^{(\lambda-3\norm{f}_\infty) T}\sigma^2\int_{a^*_{a,m}}^{T-a-\Delta_{\text{abs}}}\int_0^\infty u_0(a_0, e^{\lambda(T-a-a')}\gamma^{-1}(e^{\lambda a'}\gamma^{-1}(e^{\lambda a}m)))da_0da'.
\end{equation}
We have obtained that on $[0,\bar{a}]\times[\underline{m},\overline{m}]$, the density $ (S^{\tilde{x}}_{T}u_0)$ is lower bounded by a constant depending on a Lebesgue integral on $\R_+^2$ involving $u_0$.

\textbf{Step~2:}

Now, we want express the lower bound Eq.~\eqref{eq:lower_bound_rho_T} in terms of $\int_{\R_+\times]0,R]} u_0\,dadm$ by a change of variable. Let us define the function $\psi^T_{a,m}$:
\begin{equation*}
    \psi^T_{a,m}:[a^*_{a,m},T-a-\Delta_{\text{abs}}] \to \R_+, \quad a' \mapsto e^{\lambda(T-a-a')}\gamma^{-1}(e^{\lambda a'}\gamma^{-1}(e^{\lambda a}m)).
\end{equation*}
We verify that $(\psi^T_{a,m})'>0$: 

For all $a'\in[a^*_{a,m},T-a]$,
\begin{align}\label{eq:psi_derivative}
    (\psi^T_{a,m})'(a') = \lambda e^{\lambda(T-a-a')}\bigg\{(\gamma^{-1})'(e^{\lambda a'}\gamma^{-1}(e^{\lambda a}m))e^{\lambda a'}\gamma^{-1}(e^{\lambda a}m)- \gamma^{-1}(e^{\lambda a'}\gamma^{-1}(e^{\lambda a}m))\bigg\}.
\end{align}
As $\Gamma > 0$ and $\gamma'\leq 1$ (Assumption~\ref{assumption:long-time}), we have
\begin{align*}
    (\psi^T_{a,m})'(a') &> \lambda e^{\lambda(T-a-a')}\bigg\{(\gamma^{-1})'(e^{\lambda a'}\gamma^{-1}(e^{\lambda a}m))e^{\lambda a'}\gamma^{-1}(e^{\lambda a}m)- e^{\lambda a'}\gamma^{-1}(e^{\lambda a}m)\bigg\} \\
    &=\lambda e^{\lambda(T-a)}\gamma^{-1}(e^{\lambda a}m)\bigg\{\underbrace{(\gamma^{-1})'(e^{\lambda a'}\gamma^{-1}(e^{\lambda a}m))}_{\geq 1}- 1\bigg\}\geq 0.
\end{align*}
Therefore, $\psi^T_{a,m}$ is a strictly increasing $\C^1$-diffeomorphism from $[a^*_{a,m},T-a-\Delta_{\text{abs}}]$ to $[\psi^T_{a,m}(a^*_{a,m}), \psi^T_{a,m}(T-a-\Delta_{\text{abs}})]$. We can now rewrite Eq.~\eqref{eq:lower_bound_rho_T}:
\begin{align*}
    (S^{\tilde{x}}_{T}u_0)(a,m)&\geq e^{(\lambda-3\norm{f}_\infty) T}\sigma^2\int_{a^*_{a,m}}^{T-a-\Delta_{\text{abs}}}\int_0^\infty u_0(a_0, \psi^T_{a,m}(a'))da_0da'\\
    &=e^{(\lambda-3\norm{f}_\infty) T}\sigma^2\int_{\psi^T_{a,m}(a^*_{a,m})}^{\psi^T_{a,m}(T-a-\Delta_{\text{abs}})}\int_0^\infty u_0(a_0, m_0)\left|((\psi^T_{a,m})^{-1})'(m_0)\right|da_0dm_0.
\end{align*}
Going back to Eq.~\eqref{eq:psi_derivative}, and using the fact that there exists $C_\gamma$ such that $C_\gamma\leq \gamma'\leq 1$ (Assumption~\ref{assumption:long-time}), we have, for all $a'\in[a^*_{a,m},T-a-\Delta_{\text{abs}}]$,
\begin{align*}
    (\psi^T_{a,m})'(a') &\leq \lambda e^{\lambda(T-a-a')}C_\gamma^{-1}e^{\lambda a'}\gamma^{-1}(e^{\lambda a}m) \leq \lambda e^{\lambda T}C_\gamma^{-1}\overline{m}.
\end{align*}
Hence, 
\begin{equation*}
    (S^{\tilde{x}}_{T}u_0)(a,m)\geq \frac{e^{-3\norm{f}_\infty T}\sigma^2C_\gamma}{\lambda \overline{m}}\int_{\psi^T_{a,m}(a^*_{a,m})}^{\psi^T_{a,m}(T-a-\Delta_{\text{abs}})}\int_0^\infty u_0(a_0, m_0)da_0dm_0.
\end{equation*}
In addition, by our choice of $\{T, \bar{a}, \underline{m}, \overline{m}\}$, we have
\begin{align*}
    &\psi^T_{a,m}(a^*_{a,m}) = 0, \\
    &\psi^T_{a,m}(T-a-\Delta_{\text{abs}}) = e^{\lambda \Delta_{\text{abs}}}\gamma^{-1}(e^{\lambda(T-a-\Delta_{\text{abs}})}\gamma^{-1}(e^{\lambda a}m)) > R.
\end{align*}
Therefore,
\begin{equation*}
    (S^{\tilde{x}}_{T}u_0)(a,m)\geq \frac{e^{-3\norm{f}_\infty T}\sigma^2C_\gamma}{\lambda \overline{m}}\int_0^R\int_0^\infty u_0(a_0, m_0)da_0dm_0.
\end{equation*}
Since we have supposed that $(a,m)\in[0,\bar{a}]\times[\underline{m},\overline{m}]$, this concludes the proof.
\end{proof}

With the Lyapunov condition~\eqref{eq:lyapunov} and the Doeblin minoration condition~\eqref{eq:doeblin}, we can apply a version of Harris' theorem:
\begin{theorem}[Harris]\label{theorem:harris}
Grant Assumptions~\ref{assumption:functions} -- \ref{assumption:long-time}. For all $\tilde{x}\in\R$, there exists a unique $\rho^{\tilde{x}}_{\infty}\in L^1_+(w)$ with $\norm{\rho^{\tilde{x}}_{\infty}}_{L^1}=1$ such that $S^{\tilde{x}}_t\rho^{\tilde{x}}_{\infty}=\rho^{\tilde{x}}_{\infty}$, for all $t\geq0$, and there exists $K\geq 1$ and $\mathfrak{a}>0$ such that for all initial data $u_0\in L^1_+(w)$ with $\norm{u_0}_{L^1}=1$,
\begin{equation}
    \norm{S_t^{\tilde{x}}u_0 - \rho^{\tilde{x}}_{\infty}}_{L^1(w)} \leq Ke^{-\mathfrak{a}t}\norm{u_0-\rho^{\tilde{x}}_{\infty}}_{L^1(w)}, \quad \forall t\geq 0.
\end{equation}
Furthermore, by Lemma~\ref{lemma:bounds_l1w}, we have that $\norm{\rho^{\tilde{x}}_{\infty}}_{L^1(w)}\leq\frac{b}{\alpha}$, where the constants $\alpha$ and $b$ are taken from the Lyapunov condition~\eqref{eq:lyapunov}.
\end{theorem}
\begin{proof}
This is a classic result which proof can be found in the work of Hairer and Mattingly \cite{HaiMat11}.
\end{proof}

If, in addition, Assumption~\ref{assumption:compact} holds, we can simply apply Doeblin's theorem:
\begin{theorem}[Doeblin]\label{theorem:doeblin}
Grant Assumptions~\ref{assumption:functions}, \ref{assumption:long-time} and \ref{assumption:compact}. For all $\tilde{x}\in\R$, there exists a unique $\rho^{\tilde{x}}_{\infty}\in L^1_+$ with $\norm{\rho^{\tilde{x}}_{\infty}}_{L^1}=1$ such that $S^{\tilde{x}}_t\rho^{\tilde{x}}_{\infty}=\rho^{\tilde{x}}_{\infty}$, for all $t\geq0$, and there exists $K\geq 1$ and $\mathfrak{a}>0$ such that for all initial data $u_0\in L^1_+$ with $\norm{u_0}_{L^1}=1$,
\begin{equation}
    \norm{S_t^{\tilde{x}}u_0 - \rho^{\tilde{x}}_{\infty}}_{L^1} \leq Ke^{-\mathfrak{a}t}\norm{u_0-\rho^{\tilde{x}}_{\infty}}_{L^1}, \quad \forall t\geq 0.
\end{equation}
\end{theorem}
\begin{proof} See, for example, Theorem 2.3 in \cite{CanYol19}.
\end{proof}

We say that $\rho^{\tilde{x}}_{\infty}$ is the \textit{invariant probability measure} of the semigroup $(S^{\tilde{x}}_t)_{t\in\R_+}$. Note that both theorems imply the next corollary.
\begin{corollary}\label{corollary:weak_invariant}
Grant Assumptions~\ref{assumption:functions} -- \ref{assumption:long-time}. For all $\tilde{x}\in\R$, there exists a unique $\rho^{\tilde{x}}_{\infty}\in L^1_+(w)$ with $\norm{\rho^{\tilde{x}}_{\infty}}_{L^1}=1$ solving
\begin{subequations}\label{eq:weak_invariant}
\begin{align} 
    &\partial_a\rho^{\tilde{x}}_\infty(a,m) - \lambda\partial_m(m\rho^{\tilde{x}}_\infty(a,m)) = - f(a,m,\tilde{x})\rho^{\tilde{x}}_\infty(a,m), \\
    &\rho^{\tilde{x}}_\infty(0,m) = \mathbbm{1}_{m>\gamma(0)}\left|(\gamma^{-1})'(m)\right|\int_0^\infty f(a,\gamma^{-1}(m),\tilde{x})\rho^{\tilde{x}}_\infty(a,\gamma^{-1}(m))da, \label{eq:weak_invariant_rho}
\end{align}
\end{subequations}
in the weak sense. Furthermore, we have that $\rho^{\tilde{x}}_\infty\in\mathcal{C}(\R_+,L^1_+(\R_+^*))\cap L^\infty(\R_+,L^1_+(\R_+^*))$.
\end{corollary}

\section{Stationary solutions for arbitrary connectivity strength}\label{sec:stationary}

In this section, we study the stationary solutions to~\eqref{eq:PDE}, namely the solution to
\begin{subequations}\label{eq:stationary}
\begin{align} 
    &\partial_a\rho_\infty(a,m) - \lambda\partial_m(m\rho_\infty(a,m)) = - f(a,m,\varepsilon x_\infty)\rho_\infty(a,m),  \\
    &\rho_\infty(0,m) = \mathbbm{1}_{m>\gamma(0)}\left|(\gamma^{-1})'(m)\right|\int_0^\infty f(a,\gamma^{-1}(m),\varepsilon x_\infty)\rho_\infty(a,\gamma^{-1}(m))da,  \\
    &x_\infty = \int_{0}^\infty\int_0^\infty \bar{h}(a,m)f(a,m,\varepsilon x_\infty)\rho_\infty(a,m)dadm.
\end{align}
\end{subequations}

\begin{definition}
$(\rho_\infty,x_\infty)\in L^1_+(w)\cap\mathcal{C}(\R_+,L^1_+(\R_+^*))\cap L^\infty(\R_+,L^1_+(\R_+^*))\times\R_+$ is a \textit{stationary solution} to~\eqref{eq:PDE} if $\norm{\rho_\infty}_{L^1}=1$ and if it solves Eq.~\eqref{eq:stationary} in the weak sense.
\end{definition}

\subsection{Existence and uniqueness using the Doeblin-Harris method.}

We present two Lipschitz continuity results, which will allow us to prove the existence (and the uniqueness when $\varepsilon$ is small) of stationary solutions.
The following lemma plays the same role as Theorem 4.5 in \cite{CanYol19}:

\begin{lemma}[Lipschitz continuity at finite $T$]\label{lemma:continuity_finite_time}
Grant Assumptions~\ref{assumption:functions} -- \ref{assumption:long-time}. For all initial data $u_0\in L^1_+(w)$ and for all $T>0$, there exists a constant $C_{T,\norm{u_0}_{L^1(w)}}>0$ such that
\begin{equation}
    \forall \widetilde{x_1},\widetilde{x_2}\in \R, \qquad \norm{S_T^{\widetilde{x_1}}u_0 - S_T^{\widetilde{x_2}}u_0}_{L^1(w)} \leq C_{T,\norm{u_0}_{L^1(w)}}|\widetilde{x_1} - \widetilde{x_2}|.
\end{equation}
\end{lemma}
\begin{proof}
For all $t>0$,
\begin{align*}
    &\norm{S_t^{\widetilde{x_1}}u_0 - S_t^{\widetilde{x_2}}u_0}_{L^1(w)} \\
    &\qquad= \int_0^\infty \int_t^\infty \bigg|u_0(a-t,e^{\lambda t}m)\exp\left(\lambda t - \int_0^t f(a-t+s,e^{\lambda(t-s)}m,\widetilde{x_1})ds\right) \\
    &\qquad\qquad\qquad\qquad- u_0(a-t,e^{\lambda t}m)\exp\left(\lambda t - \int_0^t f(a-t+s,e^{\lambda(t-s)}m,\widetilde{x_2})ds\right)\bigg|w(a,m)dadm \\
    &\quad\qquad+\int_0^\infty \int_0^t \bigg|\rho^{\widetilde{x_1}}_{t-a}(0,e^{\lambda a}m)\exp\left(\lambda a - \int_{t-a}^t f(a-t+s,e^{\lambda(t-s)}m,\widetilde{x_1})ds\right) \\
    &\qquad\qquad\qquad\qquad-\rho^{\widetilde{x_2}}_{t-a}(0,e^{\lambda a}m)\exp\left(\lambda a - \int_{t-a}^t f(a-t+s,e^{\lambda(t-s)}m,\widetilde{x_2})ds\right)\bigg|w(a,m)dadm \\
    &\qquad=:Q_1 + Q_2.
\end{align*}
\begin{align*}
    Q_1 &= \int_0^\infty \int_0^\infty u_0(a,m)\left| \exp\left(-\int_0^t f(a+s,e^{-\lambda s}m,\widetilde{x_1})ds\right) - \exp\left(-\int_0^t f(a+s,e^{-\lambda s}m,\widetilde{x_2})ds\right)\right|\\
    &\qquad\qquad\qquad\qquad\qquad\qquad\qquad\qquad\qquad\qquad\qquad\qquad\qquad\qquad\qquad\qquad w(a+t,e^{-\lambda t}m)dadm \\
    &\leq \int_0^\infty \int_0^\infty u_0(a,m)\left(\int_0^t\left| f(a+s,e^{-\lambda s}m,\widetilde{x_1})- f(a+s,e^{-\lambda s}m,\widetilde{x_2})\right|ds\right) w(a+t,e^{-\lambda t}m)dadm \\
    &\leq t L_f |\widetilde{x_1}-\widetilde{x_2}|\int_0^\infty \int_0^\infty u_0(a,m) w(a+t,e^{-\lambda t}m)dadm \\
    &\leq t L_f\norm{u_0}_{L^1(w)}|\widetilde{x_1}-\widetilde{x_2}|,
\end{align*}
where in the last inequality we used
\begin{equation}\label{eq:bound_w_1}
w(a+t,e^{-\lambda t}m) \leq w(a,m),\quad \forall a\geq0, m\geq 0.
\end{equation}
\begin{align*}
    Q_2 &= \int_0^\infty \int_0^t \bigg|\rho^{\widetilde{x_1}}_{t-a}(0,m)\exp\left(- \int_{t-a}^t f(a-t+s,e^{\lambda(t-s-a)}m,\widetilde{x_1})ds\right) \\
    &\qquad\qquad\qquad-\rho^{\widetilde{x_2}}_{t-a}(0,m)\exp\left( - \int_{t-a}^t f(a-t+s,e^{\lambda(t-s-a)}m,\widetilde{x_2})ds\right)\bigg|w(a,e^{-\lambda a}m)dadm. \\
    %&= \int_0^\infty \int_0^t \bigg|\rho^{\widetilde{x_1}}_{t-a}(0,m)\exp\left(- \int_0^a f(u,e^{-\lambda u}m,\widetilde{x_1})du\right) \\
    %&\qquad\qquad\qquad-\rho^{\widetilde{x_2}}_{t-a}(0,m)\exp\left( - \int_0^a f(u,e^{-\lambda u}m,\widetilde{x_2})du\right)\bigg|w(a,e^{-\lambda a}m)dadm
\end{align*}
By changes of variables,
\begin{align*}    
    Q_2 &= \int_0^\infty \int_0^t \bigg|\rho^{\widetilde{x_1}}_{s}(0,m)\exp\left(- \int_0^{t-s} f(u,e^{-\lambda u}m,\widetilde{x_1})du\right) \\
    &\qquad\qquad\qquad-\rho^{\widetilde{x_2}}_{s}(0,m)\exp\left( - \int_0^{t-s} f(u,e^{-\lambda u}m,\widetilde{x_2})du\right)\bigg|w(t-s,e^{-\lambda (t-s)}m)dsdm \\
    &\leq \int_0^\infty \int_0^t \rho^{\widetilde{x_1}}_{s}(0,m)\bigg|\exp\left(- \int_0^{t-s} f(u,e^{-\lambda u}m,\widetilde{x_1})du\right) - \exp\left(- \int_0^{t-s} f(u,e^{-\lambda u}m,\widetilde{x_2})du\right) \bigg|\\
    &\qquad\qquad\qquad\qquad\qquad\qquad\qquad\qquad\qquad\qquad\qquad\qquad\qquad\qquad\qquad w(t-s,e^{-\lambda(t-s)}m)dsdm \\
    &\qquad+\int_0^\infty \int_0^t \left|\rho^{\widetilde{x_1}}_{s}(0,m) - \rho^{\widetilde{x_2}}_{s}(0,m)\right|w(t-s,e^{-\lambda(t-s)}m)dsdm \\
    &=: Q_{2,1} + Q_{2,2}
\end{align*}
\begin{align*}
    Q_{2,1} &\leq t\norm{f}_\infty L_f|\widetilde{x_1}-\widetilde{x_2}|  \int_0^t\int_0^\infty \int_0^\infty \left|(\gamma^{-1})' (m)\right|\rho^{\widetilde{x_1}}_{s}(a,\gamma^{-1}(m))w(t,m)dadmds \\
    &\leq t\norm{f}_\infty L_f|\widetilde{x_1}-\widetilde{x_2}|  \int_0^t\int_0^\infty \int_0^\infty \rho^{\widetilde{x_1}}_{s}(a,m)w(t,m+\norm{\Gamma}_\infty)dadmds \\
    &\leq t(1+\norm{\Gamma}_\infty)\norm{f}_\infty L_f|\widetilde{x_1}-\widetilde{x_2}|  \int_0^t \norm{\rho^{\widetilde{x_1}}_{s}}_{L^1(w)}ds,
\end{align*}
where in the last inequality we used 
\begin{equation}\label{eq:bound_w_2}
    w(t,m+\norm{\Gamma}_\infty)=1+m+\norm{\Gamma}_\infty \leq (1+\norm{\Gamma}_\infty)w(a,m), \quad \forall a\geq 0, m\geq 0.
\end{equation}
By Lemma~\ref{lemma:bounds_l1w},
\begin{align*}
    Q_{2,1} &\leq t^2(1+\norm{\Gamma}_\infty)\norm{f}_\infty L_f  \left(\norm{u_0}_{L^1(w)}+\frac{b}{\alpha}\right)|\widetilde{x_1}-\widetilde{x_2}|.
\end{align*}
\begin{align*}
    Q_{2,2} &\leq \norm{f}_\infty \int_0^t \int_0^\infty \int_0^\infty \left|(\gamma^{-1})'(m)\right|\left|\rho^{\widetilde{x_1}}_{s}(a,\gamma^{-1}(m))-\rho^{\widetilde{x_2}}_{s}(a,\gamma^{-1}(m))\right| w(t,m)dadmds \\
    &\leq \norm{f}_\infty \int_0^t \int_0^\infty \int_0^\infty \left|\rho^{\widetilde{x_1}}_{s}(a,m)-\rho^{\widetilde{x_2}}_{s}(a,m)\right| w(t,m+\norm{\Gamma}_\infty)dadmds \\
    &\leq (1+\norm{\Gamma}_\infty)\norm{f}_\infty \int_0^t \norm{S^{\widetilde{x_1}}_s u_0 - S^{\widetilde{x_2}}_s u_0}_{L^1(w)}ds,
\end{align*}
where again, in the last inequality, we used Eq.~\eqref{eq:bound_w_2}. Fix $T>0$. Gathering the bounds for $Q_1$, $Q_{2,1}$ and $Q_{2,2}$ we see that there exists constants $C>0$ and $C'_{T,\norm{u_0}_{L^1(w)}}>0$ such that, for all $t\in[0,T]$,
\begin{equation*}
    \norm{S_t^{\widetilde{x_1}}u_0 - S_t^{\widetilde{x_2}}u_0}_{L^1(w)} \leq C \int_0^t \norm{S_s^{\widetilde{x_1}}u_0 - S_s^{\widetilde{x_2}}u_0}_{L^1(w)}ds + tC'_{T,\norm{u_0}_{L^1(w)}}|\widetilde{x_1}-\widetilde{x_2}|.
\end{equation*}
By Grönwall's lemma, for all $t\in[0,T]$,
\begin{equation} \label{eq:gronwall_lipschitz}
    \norm{S_t^{\widetilde{x_1}}u_0 - S_t^{\widetilde{x_2}}u_0}_{L^1(w)} \leq \frac{C'_{T,\norm{u_0}_{L^1(w)}}|\widetilde{x_1}-\widetilde{x_2}|}{C}\left(\exp(C t)-1\right).
\end{equation}
Since Eq.~\eqref{eq:gronwall_lipschitz} holds for all $t\in[0,T]$, this achieves the proof.
\end{proof}

\begin{lemma}[Lipschitz continuity at $T=\infty$]\label{lemma:continuity}
Grant Assumptions~\ref{assumption:functions} -- \ref{assumption:long-time}. Writing $\rho^{\tilde{x}}_\infty\in L^1_+(w)$ the invariant probability measure given by Theorem~\ref{theorem:harris} for any $\tilde{x}\in\R$, the function
\begin{align*}
    \Upsilon : \R_+ \to \R_+, \qquad \Upsilon(x) = \int_0^\infty \int_0^\infty \bar{h}(a,m)f(a,m,\varepsilon x)\rho^{\varepsilon x}_\infty(a,m)dadm
\end{align*}
is Lipschitz and there exists $C>0$ such that
\begin{equation*}
    \forall x_1,x_2 \in \R_+, \qquad |\Upsilon(x_1) - \Upsilon(x_2)|\leq |\varepsilon|C|x_1-x_2|.
\end{equation*}
\end{lemma}

\begin{proof} Since $f$ is Lipschitz in $x$, we have, for any $x_1,x_2\in\R_+$,
\begin{align*}
    |\Upsilon(x_1) - \Upsilon(x_2)| &\leq \norm{\bar{h}}_\infty\left\{\norm{f}_\infty \norm{\rho^{\varepsilon x_1}_\infty - \rho^{\varepsilon x_2}_\infty}_{L^1} + L_f|\varepsilon||x_1 - x_2|\right\} \\
    &\leq \norm{\bar{h}}_\infty\left\{\norm{f}_\infty \norm{\rho^{\varepsilon x_1}_\infty - \rho^{\varepsilon x_2}_\infty}_{L^1(w)} + L_f|\varepsilon||x_1 - x_2|\right\},
\end{align*}
from where we only need to bound the first term on the right hand side. We can use Theorem~\ref{theorem:harris} and Lemma~\ref{lemma:continuity_finite_time}: for any $T\in \R_+$,
\begin{align*}
    \norm{\rho^{\varepsilon x_1}_\infty - \rho^{\varepsilon x_2}_\infty}_{L^1(w)} &=\|S_T^{\varepsilon x_1}\rho_{\infty}^{\varepsilon x_1}-S_T^{\varepsilon x_1}\rho_{\infty}^{\varepsilon x_2}+S_T^{\varepsilon x_1}\rho_{\infty}^{\varepsilon x_2}-S_T^{\varepsilon x_2}\rho_{\infty}^{\varepsilon x_2}\|_{L^1(w)} \\&\leq Ke^{-\mathfrak{a} T}\|\rho_\infty^{\varepsilon x_1}-\rho_\infty^{\varepsilon x_2}\|_{L^1(w)}+C_T |\varepsilon||x_1-x_2|,
\end{align*}
where $K$ and $\mathfrak{a}$ are the exponential stability constants of Theorem~\ref{theorem:harris}.  Choosing $T$ such that $Ke^{-\mathfrak{a}T} = 1/2$, we get
\begin{equation*}
    \norm{\rho^{\varepsilon x_1}_\infty - \rho^{\varepsilon x_2}_\infty}_{L^1(w)} \leq 2C_T |\varepsilon||x_1-x_2|.
\end{equation*}
Gathering the bounds concludes the proof.
\end{proof}

\begin{theorem}[Stationary solutions]
Grant Assumptions~\ref{assumption:functions} -- \ref{assumption:long-time}. We have
\begin{enumerate}[label=(\roman*)]
    \item There exists a stationary solution to~\eqref{eq:PDE}.
    \item There exists $\varepsilon^*>0$ such that for all $\varepsilon\in]-\varepsilon^*,+\varepsilon^*[$, the stationary solution to~\eqref{eq:PDE} is unique.
\end{enumerate}
\end{theorem}
\begin{proof}
For all $\tilde{x}\in\R$, let us write $\rho^{\tilde{x}}_\infty\in L^1_+(w)$ the unique invariant measure given by Theorem~\ref{theorem:harris} and let us also take the function $\Upsilon$ from Lemma~\ref{lemma:continuity}. By Corollary~\ref{corollary:weak_invariant}, $(\rho_\infty,x_\infty)\in L^1_+(w)\cap\mathcal{C}(\R_+,L^1_+(\R_+^*))\cap L^\infty(\R_+,L^1_+(\R_+^*))\times\R_+$ is a weak solution to Eq.~\eqref{eq:stationary} if and only if $\rho_\infty = \rho^{\varepsilon x_\infty}_\infty$ and $x_\infty$ is a fixed-point of $\Upsilon$. Hence, the study of the existence and the uniqueness of stationary solutions is reduced to the study of the existence and the uniqueness of the fixed-point of $\Upsilon$.

Since for all $x\in\R_+$, $\norm{\rho^{\varepsilon x}_\infty}_{L^1}=1$, we have that for all $x\in\R_+$, $\Upsilon(x)\leq \norm{\bar{h}}_\infty\norm{f}_\infty$. Therefore, the set $[0,\norm{\bar{h}}_\infty\norm{f}_\infty]$ (which is compact and convex) is stable by $\Upsilon$. { Then, the continuity of $\Upsilon$ guarantees the existence of a fixed-point, which proves $(i)$.} 

To obtain $(ii)$, we observe that the Lipschitz constant of $\Upsilon$ is $|\varepsilon|C$: if we take $|\varepsilon|< \varepsilon^*:= C^{-1}$, $\Upsilon$ is a contraction and we can apply Banach's fixed-point theorem to conclude.
\end{proof}

\subsection{Alternative proof for the existence using Schauder's fixed-point theorem}

We include here an alternative proof for the existence of a stationary solution, which is interesting for two reasons: on the one hand, it does not rely on the Harris-Doeblin method, and on the other hand, it provides some estimates on the stationary solutions.

For any $(\tilde{u},\tilde{x})\in L^1_+(]\gamma(0),+\infty[)\times\R$, consider the transport equation 
\begin{align*}
    &\partial_a\varrho(a,m) - \lambda\partial_m(m\varrho(a,m)) = - f(a,m,\tilde{x})\varrho(a,m),\\
    &\varrho(0,m) = \tilde{u}(m).
\end{align*}
It has a unique weak solution $\rho_\infty^{\tilde{u},\tilde{x}} \in \mathcal{C}(\R_+,L^1_+(\R_+^*))\cap L^\infty(\R_+,L^1_+(\R_+^*))$ given by the method of characteristics: for all $(a,m)\in\R_+\times\R_+^*$,
\begin{equation} \label{eq:stationary_sol}
    \rho_\infty^{\tilde{u},\tilde{x}}(a,m) = \tilde{u}(e^{\lambda a}m)\exp\left(\lambda a - \int_0^a f(s,e^{\lambda(a-s)}m,\tilde{x})ds\right).
\end{equation}

We can now define the operator $\Phi := (\Phi_1, \Phi_2)$ on $L^1_+(]\gamma(0),+\infty[)\times\R$ where, for all $(\tilde{u},\tilde{x})\in L^1_+(]\gamma(0),+\infty[)\times\R$,
\begin{subequations}\label{eq:def_phi}
\begin{align}
    \Phi_1(\tilde{u},\tilde{x})(m) &:= \mathbbm{1}_{m>\gamma(0)}\left|(\gamma^{-1})'(m)\right|\int_0^\infty f(a,\gamma^{-1}(m),\tilde{x})\rho_\infty^{(\tilde{u},\tilde{x})}(a,\gamma^{-1}(m))da, \\
    \Phi_2(\tilde{u},\tilde{x}) &:= \int_{0}^\infty\int_0^\infty \bar{h}(a,m)f(a,m,\tilde{x})\rho_\infty^{\tilde{u},\tilde{x}}(a,m)dadm. \label{eq:phi2}
\end{align}
\end{subequations}
$(\rho_\infty,x_\infty)$ is a stationary solution if and only if it is a fixed-point of $\Phi$. Whence, we get the \textit{a priori} estimates:

\begin{lemma}\label{lemma:a_priori_stationary}
Grant Assumptions~\ref{assumption:functions} and \ref{assumption:long-time}. There exists $\theta \in ]0,1[$ such that for all $(\tilde{u},\tilde{x})\in L^1_+(]\gamma(0),+\infty[)\times\R$, 
\begin{enumerate}[label=(\roman*)]
\item $\norm{\Phi_1(\tilde{u},\tilde{x})}_{L^1} = \norm{\tilde{u}}_{L^1}$.
\item For all $m\in\R_+^*$, $|\Phi_1(\tilde{u},\tilde{x})(m)|\leq \mathbbm{1}_{m>\gamma(0)} \frac{\norm{f}_\infty}{\lambda \gamma^{-1}(m)}\norm{\tilde{u}}_{L^1}$.
\item  
\begin{equation*}
    \int_0^\infty \Phi_1(\tilde{u},\tilde{x})(m)mdm \leq \max\left(\int_0^\infty \tilde{u}(m)mdm, \frac{\gamma(0)}{1-\theta}\norm{\tilde{u}}_{L^1}\right).
\end{equation*}
\item For all $\beta \in ]0,\frac{\min(f)}{\lambda}[$,
\begin{equation*}
    \int_{\gamma(0)}^\infty\frac{\Phi_1(\tilde{u},\tilde{x})(m)}{\gamma^{-1}(m)^\beta}dm \leq \frac{\norm{f}_\infty}{\lambda\gamma(0)^\beta}\left(\frac{\min(f)}{\lambda} - \beta\right)\norm{\tilde{u}}_{L^1}.
\end{equation*}
\item $\Phi_2(\tilde{u},\tilde{x})\leq \norm{\bar{h}}_\infty\norm{\tilde{u}}_{L^1}$.
\end{enumerate}
\end{lemma}

\begin{proof}
$(i)$ By changes of variables on $m$,
\begin{align*}
    \norm{\Phi_1(\tilde{u},\tilde{x})}_{L^1} &= \int_0^\infty\int_0^\infty f(a,m,\tilde{x})\tilde{u}(e^{\lambda a}m)\exp\left(\lambda a - \int_0^a f(s,e^{\lambda(a-s)}m,\tilde{x})ds\right)dadm \\
    &=\int_0^\infty \tilde{u}(m)\underbrace{\int_0^\infty f(a,e^{-\lambda a}m,\tilde{x})\exp\left(-\int_0^a f(s,e^{-\lambda s}m,\tilde{x})ds\right)da}_{=1\quad { (\text{by Assumption~\ref{assumption:long-time}~$(i)$)}}}dm.
\end{align*}

$(ii)$
\begin{align*}
    |\Phi_1(\tilde{u},\tilde{x})(m)| &\leq \mathbbm{1}_{m>\gamma(0)}\norm{f}_\infty \int_0^\infty \tilde{u}(e^{\lambda a}\gamma^{-1}(m))\exp(\lambda a)da \\
    &=\mathbbm{1}_{m>\gamma(0)}\frac{\norm{f}_\infty}{\lambda \gamma^{-1}(m)}\int_0^\infty \tilde{u}(e^{\lambda a}\gamma^{-1}(m))\gamma^{-1}(m)\lambda\exp(\lambda a)da \\
    &=\mathbbm{1}_{m>\gamma(0)}\frac{\norm{f}_\infty}{\lambda \gamma^{-1}(m)}\underbrace{\int_{\gamma^{-1}(m)}^\infty \tilde{u}(y)dy}_{\leq \norm{\tilde{u}}_{L^1}},
\end{align*}
where for the last equality we used the change of variable $y=e^{\lambda a}\gamma^{-1}(m)$.

$(iii)$ Performing the same change of variable as for $(i)$ and using the fact that $\gamma(m) \leq \gamma(0) + m$, $\forall m\in \R_+$ (since $\gamma' \leq 1$), we have
\begin{align*}
    \int_0^\infty \Phi_1&(\tilde{u},\tilde{x})(m)mdm \\
    &= \int_0^\infty \tilde{u}(m)\int_0^\infty \gamma(e^{-\lambda a}m)f(a,e^{-\lambda a}m,\tilde{x})\exp\left(-\int_0^a f(s,e^{-\lambda s}m,\tilde{x})ds\right)dadm \\
    &\leq \int_0^\infty \tilde{u}(m)m \underbrace{\int_0^\infty e^{-\lambda a}f(a,e^{-\lambda a}m,\tilde{x})\exp\left(-\int_0^a f(s,e^{-\lambda s}m,\tilde{x})ds\right)da}_{=:\vartheta(m)}dm + \gamma(0)\norm{\tilde{u}}_{L^1}.
\end{align*}

There exists $\theta\in]0,1[$ such that for all $m\in\R_+^*$, $\vartheta(m)<1$: \\
Fix $\epsilon>0$.
\begin{align*}
    \vartheta(m) &\leq \int_0^\epsilon f(a,e^{-\lambda a}m,\tilde{x})\exp\left(-\int_0^a f(s,e^{-\lambda s}m,\tilde{x})ds\right)da \\
    &\quad + \int_{\epsilon}^\infty e^{-\lambda \epsilon}f(a,e^{-\lambda a}m,\tilde{x})\exp\left(-\int_0^a f(s,e^{-\lambda s}m,\tilde{x})ds\right)da\\
    &= 1-(1-e^{-\lambda \epsilon})\int_{\epsilon}^\infty f(a,e^{-\lambda a}m,\tilde{x})\exp\left(-\int_0^a f(s,e^{-\lambda s}m,\tilde{x})ds\right)da\\
    &= 1-(1-e^{-\lambda \epsilon})\exp\left(-\int_0^\epsilon f(s,e^{-\lambda s}m,\tilde{x})ds\right) \\
    &\leq 1-(1-e^{-\lambda \epsilon})\exp(-\norm{f}_\infty \epsilon)  =: \theta <1.
\end{align*}
Whence,
\begin{equation*}
    \int_0^\infty \Phi_1(\tilde{u},\tilde{x})(m)mdm \leq \theta \int_0^\infty \tilde{u}(m)mdm + \gamma(0)\norm{\tilde{u}}_{L^1}.
\end{equation*}
To see that 
\begin{equation*}
    \int_0^\infty \Phi_1(\tilde{u},\tilde{x})(m)mdm \leq \max\left(\int_0^\infty\tilde{u}(m)mdm, \frac{\gamma(0)}{1-\theta}\norm{\tilde{u}}_{L^1}\right),
\end{equation*}
we can distinguish three cases: if $\int_0^\infty\tilde{u}(m)mdm = \infty$,  the inequality is trivial; if $\frac{\gamma(0)}{1-\theta}\norm{\tilde{u}}_{L^1}\leq \int_0^\infty\tilde{u}(m)mdm < +\infty$, then
\begin{align*}
    \int_0^\infty \Phi_1(\tilde{u},\tilde{x})(m)mdm &\leq \int_0^\infty \tilde{u}(m)mdm - (1-\theta)\int_0^\infty \tilde{u}(m)mdm + \gamma(0)\norm{\tilde{u}}_{L^1}\\
    &\leq \int_0^\infty \tilde{u}(m)mdm;
\end{align*}
and finally if $\int_0^\infty\tilde{u}(m)mdm<\frac{\gamma(0)}{1-\theta}\norm{\tilde{u}}_{L^1}$, then
\begin{equation*}
    \int_0^\infty \Phi_1(\tilde{u},\tilde{x})(m)mdm \leq \theta \frac{\gamma(0)}{1-\theta}\norm{\tilde{u}}_{L^1} + \gamma(0)\norm{\tilde{u}}_{L^1} = \frac{\gamma(0)}{1-\theta}\norm{\tilde{u}}_{L^1}.
\end{equation*}

$(iv)$ 
\begin{align*}
    \int_{\gamma(0)}^\infty\frac{\Phi_1(\tilde{u},\tilde{x})(m)}{\gamma^{-1}(m)^\beta}dm &= \int_0^\infty \int_0^\infty \frac{1}{m}f(a,m,\tilde{x})\tilde{u}(e^{\lambda a}m)\exp\left(\lambda a - \int_0^a f(s,e^{\lambda(a-s)}m,\tilde{x})\right)dadm \\
    &\leq\norm{f}_\infty \int_0^\infty \int_0^\infty \frac{1}{m^\beta}\tilde{u}(e^{\lambda a}m)\exp\left(\lambda a - \min(f)a\right)dadm,
\end{align*}
making the change of variable $y = e^{\lambda a}m$:
\begin{align*}
    &=\norm{f}_\infty \int_0^\infty \int_m^\infty \frac{1}{\lambda m^{1+\beta}}\tilde{u}(y)\exp\left( - \min(f)\frac{1}{\lambda}\ln\left(\frac{y}{m}\right)\right)dydm \\
    &=\frac{\norm{f}_\infty}{\lambda} \int_0^\infty \int_m^\infty m^{\min(f)/\lambda - 1 - \beta}\tilde{u}(y)y^{-\min(f)/\lambda}dydm,
\end{align*}
using Fubini's theorem and the fact that $\min(f)/\lambda - \beta > 0$:
\begin{align*}
    &= \frac{\norm{f}_\infty}{\lambda}\int_0^\infty \tilde{u}(y)y^{-\min(f)/\lambda}\underbrace{\int_0^y m^{\min(f)/\lambda - 1 - \beta}dm}_{=\frac{y^{\min(f)/\lambda - \beta}}{\min(f)/\lambda - \beta}} dy \\
    &=\frac{\norm{f}_\infty}{\lambda}\left(\frac{\min(f)}{\lambda} - \beta\right)\int_0^\infty \tilde{u}(y)y^{-\beta}dy.
\end{align*}
Finally, it is easy to check that $\int_0^\infty \tilde{u}(y)y^{-\beta}dy\leq \gamma(0)^{-\beta}\norm{\tilde{u}}_{L^1}$.

$(v)$ Use Eq.~\eqref{eq:phi2} and see the proof of $(i)$.
\end{proof}

By these estimates, we see that there exists $\beta,C_1,C_2,C_3,C_4>0$ such that the set $\mathscr{C}\times B\in L^1(]\gamma(0),+\infty[)\times\R$, where
\begin{multline*}
    \mathscr{C} := \biggr\{u\in L^1_+(]\gamma(0),+\infty[) \;\biggr|\; \norm{u}_{L^1}\leq 1;\; u\leq\frac{C_1}{\gamma^{-1}(\cdot)}\;\; a.e.; \\
    \int_0^\infty u(m)mdm \leq C_2;\; \int_{\gamma(0)}^\infty \frac{u(m)}{\gamma^{-1}(m)^\beta}dm\leq C_3\biggr\}
\end{multline*}
and $B:=[-C_4,+C_4]$, is stable by the operator $\Phi$.

In order to apply Schauder's fixed-point theorem, we will need

\begin{lemma}\label{lemma:convexity}
Grant Assumptions~\ref{assumption:functions} and \ref{assumption:long-time}. $\mathscr{C}$ is convex, closed and compact for the weak topology $\sigma(L^1,L^\infty)$.
\end{lemma}
\begin{proof}
It is easy to verify that $\mathscr{C}$ is convex. Since $\mathscr{C}$ is convex, if suffices to show that it is strongly closed to show that it is weakly closed. Let $u_n$ be a sequence of elements of $\mathscr{C}$ which converge strongly to $u\in L^1(]\gamma(0),+\infty[)$. By the strong convergence, $\norm{u}_{L^1}\leq 1$. We can extract a subsequence $u_{n_k}$ such that $u_{n_k}$ converges to $u$ a.e. Taking the pointwise limit, we have that $u\leq \frac{C_1}{\gamma^{-1}(\cdot)}$ a.e. Furthermore, by Fatou's lemma,
\begin{equation*}
    \int_{\gamma(0)}^\infty u(m)mdm \leq \liminf_{k\to +\infty}\int_{\gamma(0)}^\infty u_{n_k}(m)mdm\leq C_2
\end{equation*}
and
\begin{equation*}
    \int_{\gamma(0)}^\infty \frac{u(m)}{\gamma^{-1}(m)^\beta}dm \leq \liminf_{k\to +\infty}\int_{\gamma(0)}^\infty \frac{u_{n_k}(m)}{\gamma^{-1}(m)^\beta}dm\leq C_3.
\end{equation*}
Hence, $\mathscr{C}$ is strongly closed.

To show that $\mathscr{C}$ is weakly compact, we will show that
\begin{enumerate}
	\item[a.] $\sup_{u\in\mathscr{C}}\norm{u}_{L^1}<\infty$, 
	\item[b.]  $\forall \epsilon>0$, $\exists R>0$ such that $\int_R^\infty u(m)dm<\epsilon$ for all $u\in\mathscr{C}$,
	\item[c.] $\mathscr{C}$ is equi-integrale, i.e. $\forall \epsilon>0$, $\exists \delta>0$ such that for all Borel set $A\subset\R_+$ with $|A|\leq\delta$ and for all $u\in\mathscr{C}$, $\int_{A}u(m)dm\leq\epsilon$,
\end{enumerate}
and use Dunford-Pettis theorem. (a.) is clearly verified. (b.) is also verified since for all $R>0$, $\int_R^\infty u(m)dm\leq \frac{1}{R}\int_0^\infty u(m)mdm\leq\frac{C_2}{R}$. To show (c.), let us first observe that for all $\delta_1>0$,
\begin{equation*}
    \int_{\gamma(0)}^{\gamma(0)+ \delta_1} u(m)dm \leq \gamma^{-1}(\gamma(0) + \delta_1)^\beta \int_{\gamma(0)}^\infty \frac{u(m)}{\gamma^{-1}(m)^\beta}dm \leq \gamma^{-1}(\gamma(0)+\delta_1)^\beta C_3.
\end{equation*}
For any $\epsilon>0$, let us choose $\delta_1>0$ such that $\gamma^{-1}(\gamma(0)+\delta_1)^\beta C_3\leq\frac{\epsilon}{2}$. Then, for all Borel set $A\subset\R_+$ with $|A|\leq \delta$,
\begin{equation*}
    \int_A u(m)dm \leq \int_{\gamma(0)}^{\gamma(0)+\delta_1} u(m)dm + \int_{A\setminus[0,\gamma(0)+\delta_1]}u(m)dm \leq \frac{\epsilon}{2} + \delta\frac{C_1}{\gamma^{-1}(\gamma(0)+\delta_1)}.
\end{equation*}
Hence, we can choose $\delta = \min\left(\delta_1, \frac{\epsilon \gamma^{-1}(\gamma(0)+\delta_1)}{2 C_1}\right)$ and (c.) is verified. By the Dunford-Pettis theorem, $\mathscr{C}$ is weakly relatively compact. Finally, since $\mathscr{C}$ is weakly closed, $\mathscr{C}$ is weakly compact.  
\end{proof}

We can now give an alternative proof of the existence of stationary solutions to \eqref{eq:PDE} for arbitrary connectivity strength $\varepsilon$:
%\begin{theorem}
%Grant Assumptions~\ref{assumption:functions} and \ref{assumption:long-time}. There exists a stationary solution to Eq.~\eqref{eq:PDE}.
%\end{theorem}
\begin{proof}[{ Proof of Theorem~\ref{theorem:stationary_solutions}~$(i)$}]
We verify that the operator $\Phi$ is weakly continuous: { For any sequence $(u_n,x_n) \to (u,x)$ in $\mathscr{C}\times\R$ and for any $\varphi\in L^\infty(\R_+)$,
\begin{eqnarray*}
    \left|\int (\Phi_1(u_n,x_n)-\Phi_1(u,x))\varphi(m)dm\right|\leq Q^n_1+Q^n_2+Q^n_3,
\end{eqnarray*}
where
\begin{align*}
    Q_1^n:=&\left|\int_0^\infty\int_0^\infty (u_n(e^{\lambda a}m)-u(e^{\lambda a}m))\varphi(\gamma(m))e^{\lambda a}f(a,m,x)e^{-\int_0^a f(\tau, e^{\lambda(a-\tau)}m, x)d\tau} dadm\right|,\\
    Q_2^n:=&\norm{\varphi}_\infty\int_0^\infty\int_0^\infty u_n(e^{\lambda a}m) e^{\lambda a}|f(a,m,x)-f(a,m,x_n)|e^{-\int_0^a f(\tau, e^{\lambda(a-\tau)}m, x)d\tau}dadm,\\
    Q_3^n:=&\norm{\varphi}_\infty\int_0^\infty\int_0^\infty u_n(e^{\lambda a}m) e^{\lambda a}f(a,m,x_n)\left|e^{-\int_0^a f(\tau, e^{\lambda(a-\tau)}m, x)d\tau}-e^{-\int_0^a f(\tau, e^{\lambda(a-\tau)}m, x_n)d\tau}\right|dadm.
\end{align*}
Making the change of variable $ydy=e^{\lambda a}mdm$ in $Q_1$ we get  
\begin{equation*}  
Q_1^n=\left|\int_0^\infty (u_n(y)-u(y))\int_0^\infty \varphi(\gamma(ye^{-\lambda a}))f(a,y e^{-\lambda a},x)e^{-\int_0^a f(\tau, e^{-\lambda\tau}y, x)}dadm\right|. 
\end{equation*}  
Since $u_n$ converges to $u$ in $\sigma(L^1,L^\infty)$ and
\begin{multline*}
    \int_0^\infty \varphi(\gamma(ye^{-\lambda a}))f(a,y e^{-\lambda a},x)e^{-\int_0^a f(\tau, e^{-\lambda\tau}y, x)d\tau}da \\
    \leq \norm{\varphi}_\infty \int_0^\infty f(a,y e^{-\lambda a},x)e^{-\int_0^a f(\tau, e^{-\lambda\tau}y, x)d\tau}da=\norm{\varphi}_\infty,
\end{multline*}
$Q_1^n$ converges to 0. On the other hand, since $f$ is bounded and Lipschitz,  $Q_2^n,Q_3^n\leq\|u_n\|_{L^1}C|x_n-x|\leq C|x_n-x|$. Whence, $\Phi_1$ is a continuous operator with respect to the weak topology $\sigma(L^1,L^\infty)$.

The continuity of $\Phi_2$ is shown analogously, taking $\varphi=h$ ($h$ is a bounded).}

Since $\mathscr{C}$ is stable by $\Phi$, convex and weakly compact (Lemma~\ref{lemma:convexity}), we can apply Schauder's fixed-point theorem to obtain the existence of a fixed-point, which gives the existence of a stationary solution. 
\end{proof}

\begin{corollary} Grant Assumptions~\ref{assumption:functions} and \ref{assumption:long-time}. If $f$ is of class $\mathcal{C}^{k}$, then $u(m)$ is a function of class $\mathcal{C}^k$ for all $m>\gamma(0)$. Consequently, the stationary solutions of \eqref{eq:PDE} are of class $\mathcal{C}^k$.
\end{corollary}

\begin{proof}
If $(u,\tilde{x})$ is a fixed-point of $\Phi$, then
\begin{multline}\label{eq:u_formula}
    u(m)=\mathbbm{1}_{m>\gamma(0)}\left|(\gamma^{-1})'(m)\right|\int_0^\infty f(a,\gamma^{-1}(m),\tilde{x})u(e^{\lambda a}\gamma^{-1}(m))\\
    \exp\left(\lambda a - \int_0^a f(s,e^{\lambda(a-s)}\gamma^{-1}(m),\tilde{x})ds\right)da.
\end{multline}
Making the change of variable $y=e^{\lambda a}\gamma^{-1}(m)$ in $a$, as in the estimate $(ii)$ of Lemma \ref{lemma:a_priori_stationary}, we obtain
\begin{equation}\label{eq:u_formula_change}
    u(m)=\mathbbm{1}_{m>\gamma(0)}\frac{\left|(\gamma^{-1})'(m)\right|}{\lambda \gamma^{-1}(m)}\int_{\gamma^{-1}(m)}^\infty f(g(y,m),y,\tilde{x})u(y)\exp\left( - \int_0^{g(y,m)} f(s,e^{s}y,\tilde{x})ds\right)dy,
\end{equation}
where $g(y,m)=\ln{\frac{y}{\lambda(\gamma^{-1}(m))}}$. { We conclude with a bootstrap argument: if $u$ is $L^1$, then the right hand side of Eq.~\eqref{eq:u_formula_change} is a continuous function of $m$, meaning that $u$ is continuous. But if $u$ is continuous, then the right hand side is of class $\mathcal{C}^1$, etc.}
\end{proof}

\begin{corollary}
Grant Assumptions~\ref{assumption:functions} and \ref{assumption:long-time}. There exists a constant $C>0$, such that the stationary solution $\rho_\infty$ satisfies, 
\begin{equation}
    \rho_\infty(a,m)\leq \frac{ Ce^{-\sigma(a-\Delta_{abs})}}{m}.
\end{equation}
\end{corollary}
\begin{proof}
From the previous theorem it follows that there is $C$ such that $u(m)\leq C/m$, which, together with \eqref{eq:stationary_sol}, implies
\begin{equation}
    \rho_\infty(a,m)\leq \overline{C}\dfrac{e^{-\int_0^a f(s,e^{\lambda(a-s)}\gamma^{-1}(m),\tilde{x})ds}}{m}.
\end{equation}
The estimate follows from Assumption~\ref{assumption:long-time}~$(i)$.
\end{proof}

\subsection{Formula in the case of short-term synaptic depression}\label{sec:formula_STD}
In general, there is no explicit formula for the invariant probability measure solving Eq.~\eqref{eq:weak_invariant}. However, in the case of short-term synaptic depression Eq.~\eqref{eq:PDE_STD}, we can derive an explicit expression for the total postsynaptic potential
\begin{equation}
    X(\tilde{x}) := \int_0^\infty\hat{h}(t)\int_0^1\int_0^\infty (1-m)f(a,\tilde{x})\rho_\infty^{\tilde{x}}(a,m)dadmdt,
\end{equation}
for any $\tilde{x}\in\R$. This fact has been reported in the theoretical neuroscience literature \cite{RomAmi06}; we provide here a rigorous and analytic justification for it. 

For all $\tilde{x}\in\R$, let us introduce the quantities
\begin{align*}
    { I^{\tilde{x}}}&:=\int_0^\infty af(a,\tilde{x})\exp\left(-\int_0^a f(s,\tilde{x})ds\right)da = \int_0^\infty \exp\left(-\int_0^a f(s,\tilde{x})ds\right)da, \\
    P^{\tilde{x}}(\lambda)&:= \int_0^\infty e^{-\lambda a}f(a,\tilde{x})\exp\left(-\int_0^a f(s,\tilde{x})ds\right)da.
\end{align*}
$I^{\tilde{x}}$ can be interpreted as the mean inter-spike interval of a neuron receiving a constant input $\tilde{x}$. $P^{\tilde{x}}(\lambda)$ can be seen as the Laplace transform of the inter-spike interval distribution of that neuron, evaluated in $\lambda$.

\begin{proposition}
Grant Assumptions~\ref{assumption:functions} and \ref{assumption:long-time}. For all $\tilde{x}\in\R$,
\begin{equation*}
    X(\tilde{x}) = \int_0^\infty \hat{h}(t)dt\;\frac{1}{I^{\tilde{x}}}\left\{\frac{1-P^{\tilde{x}}(\lambda)}{1-\upsilon P^{\tilde{x}}(\lambda)}\right\}.
\end{equation*}
\end{proposition}

\begin{proof}
Using the method of characteristics { (i.e. combining Eqs.~\eqref{eq:weak_invariant_rho} and \eqref{eq:stationary_sol})}, we have
\begin{align*}
    1 &= \int_0^1 \int_0^\infty \rho_\infty^{\tilde{x}}(a,m)dadm = \int_0^1\int_0^\infty \mathbbm{1}_{e^{\lambda a}m<1}\rho_\infty^{\tilde{x}}(0,e^{\lambda a}m)\exp\left(\lambda a -\int_0^a f(s,\tilde{x})ds\right)dadm \\
    &=\int_0^1\int_0^\infty\rho_\infty^{\tilde{x}}(0,m)\exp\left(-\int_0^a f(s,\tilde{x})ds\right)dadm = I^{\tilde{x}}\int_0^1\rho_\infty^{\tilde{x}}(0,m)dm.
\end{align*}
Whence,
\begin{equation*}
    \int_0^1\int_0^\infty f(a,\tilde{x})\rho_\infty^{\tilde{x}}(a,m)dadm = \int_0^1\rho_\infty^{\tilde{x}}(0,m)dm = \frac{1}{I^{\tilde{x}}}.
\end{equation*}
On the other hand, 
\begin{align*}
    \int_0^1\int_0^\infty m f(a,\tilde{x})&\rho_\infty^{\tilde{x}}(a,m)dadm\\ &=\int_0^1\int_0^\infty \mathbbm{1}_{e^{\lambda a}m<1}mf(a,\tilde{x})\rho_\infty^{\tilde{x}}(0,e^{\lambda a}m)\exp\left(\lambda a - \int_0^a f(s,\tilde{x})ds\right)dadm \\
    &= \int_0^1\int_0^\infty e^{-\lambda a}mf(a,\tilde{x})\rho_\infty^{\tilde{x}}(0,m)\exp\left( - \int_0^a f(s,\tilde{x})ds\right)dadm \\
    &=P^{\tilde{x}}(\lambda)\int_0^1 m\rho_\infty^{\tilde{x}}(0,m)dm
\end{align*}
and
\begin{align*}
    \int_0^1 m\rho_\infty^{\tilde{x}}(0,m)dm &= \int_0^1 m \mathbbm{1}_{m>1-\upsilon}\frac{1}{\upsilon}\int_0^\infty f(a,\tilde{x})\rho_\infty^{\tilde{x}}\left(a,1-\frac{1-m}{\upsilon}\right)dadm\\
    &= \int_0^1(1-\upsilon+\upsilon m)\int_0^\infty f(a,\tilde{x})\rho_\infty^{\tilde{x}}(a,m)dadm \\
    &=\frac{1-\upsilon}{I^{\tilde{x}}} + \upsilon P^{\tilde{x}}(\lambda)\int_0^1 m\rho_\infty^{\tilde{x}}(0,m)dm.
\end{align*}
Whence,
\begin{equation*}
    \int_0^1 m\rho_\infty^{\tilde{x}}(0,m)dm = \frac{1-\upsilon}{I^{\tilde{x}}(1-\upsilon P^{\tilde{x}}(\lambda))}
\end{equation*}
and
\begin{equation*}
    \int_0^1\int_0^\infty m f(a,\tilde{x})\rho_\infty^{\tilde{x}}(a,m)dadm = \frac{P^{\tilde{x}}(\lambda)(1-\upsilon)}{I^{\tilde{x}}(1-\upsilon P^{\tilde{x}}(\lambda))}.
\end{equation*}
Finally, we have
\begin{align*}
    X(\tilde{x}) &= \int_0^\infty \hat{h}(t)dt\left\{\int_0^1\int_0^\infty f(a,\tilde{x})\rho_\infty^{\tilde{x}}(a,m)dadm - \int_0^1\int_0^\infty m f(a,\tilde{x})\rho_\infty^{\tilde{x}}(a,m)dadm\right\}\\
    &=\int_0^\infty \hat{h}(t)dt\;\frac{1}{I^{\tilde{x}}}\left\{1- \frac{P^{\tilde{x}}(\lambda)(1-\upsilon)}{(1-\upsilon P^{\tilde{x}}(\lambda))}\right\} = \int_0^\infty \hat{h}(t)dt\;\frac{1}{I^{\tilde{x}}}\left\{\frac{1-P^{\tilde{x}}(\lambda)}{1-\upsilon P^{\tilde{x}}(\lambda)}\right\}.
\end{align*}
\end{proof}

\section{Exponential stability in the weak connectivity regime}\label{sec:weak}

To study the long time behavior \eqref{eq:PDE} in the weak connectivity regime, we perturb the non-interacting case~\eqref{eq:semigroup}, taking $\tilde{x} = \varepsilon x_{\infty}$, where $x_\infty$ is given by the unique stationary solution to \eqref{eq:PDE} when $\varepsilon \in ]-\varepsilon^*,+\varepsilon^*[$ ($\varepsilon^*$ is taken from Theorem~\ref{theorem:stationary_solutions}~$(ii)$). In this section, we keep the small $\varepsilon$ fixed and we work under Assumptions~\ref{assumption:functions} --
\ref{assumption:long-time} and \ref{assumption:nonlinear}. We roughly follow the same line of argument as \cite[Sec.~5]{MisQui18}.

For convenience, we first rewrite~\eqref{eq:PDE} in a more formal and compact form:
\begin{subequations} \label{eq:PDE_compact}
\begin{align} 
    \partial_t\rho_t &= -\partial_a\rho_t + \lambda\partial_m(m\rho_t) - f(\varepsilon x_t)\rho_t + \delta_0^a(\gamma\circ\Pi)_*\left(f(\varepsilon x_t)\rho_t\right), \label{eq:PDE_compact_first}\\
    &x_t = \int_0^t \int h(t-s)f(\varepsilon x_s)\rho_s\,dadm\,ds, \\
    &\rho_0 = u_0, 
\end{align}
\end{subequations}
where $\delta_0^a$ indicates that (singular) mass enters in $a=0$\footnote{$\delta_0^a$ should not be confused with the Dirac distribution $\delta_{0=a}$. Using $\delta_{0=a}$, by integration by parts of weak solutions, Eq.~\eqref{eq:PDE_compact_first} should write
\begin{equation*}
    \partial_t\rho_t = -\partial_a\rho_t + \lambda\partial_m(m\rho_t) - f(\varepsilon x_t)\rho_t + \delta_{0=a}\left\{(\gamma\circ\Pi)_*\left(f(\varepsilon x_t)\rho_t\right) - \rho_t(0,\cdot)\right\}.
\end{equation*}}, $\Pi:(a,m)\mapsto m$ is the projection on $m$ and $_*$ denotes the pushforward measure. To write Eq.~\eqref{eq:PDE_compact} as an evolution equation, we introduce an auxiliary transport equation on $\R_+\times\R_+\times\R_+^*$
\begin{align*}
    \partial_t \zeta_t &= -\partial_s \zeta_t + \delta_0^s f(\varepsilon x_t)\rho_t, \\
    \zeta_0 &= 0,
\end{align*}
which solution is given by the method of characteristics:
\begin{equation*}
    \zeta_t(s) = \mathbbm{1}_{s\leq t} f(\varepsilon x_{t-s})\rho_{t-s}, \qquad \forall (t,s)\in\R_+^*\times\R_+.
\end{equation*}

Using the auxiliary equation, Eq.~\eqref{eq:PDE_compact} is equivalent to 
\begin{subequations}\label{eq:PDE_compact_auxiliary}
\begin{align} 
    \partial_t(\rho_t, \zeta_t) &= \left(-\partial_a\rho_t + \lambda\partial_m(m\rho_t) - f(\varepsilon x_t)\rho_t + \delta_0^a(\gamma\circ\Pi)_*\left(f(\varepsilon x_t)\rho_t\right), -\partial_s \zeta_t + \delta_0^s f(\varepsilon x_t)\rho_t \right), \\
    (\rho_0,\zeta_0) &= \left(u_0,0 \right), 
\end{align}
\end{subequations}
where $x_t := \int_0^\infty  \int h(s)\zeta_t(s)\,dadm\,ds$.

By Theorem~\ref{theorem:stationary_solutions}, for all $\varepsilon\in]-\varepsilon^*,+\varepsilon^*[$, there exists a unique stationary solution $(\rho_\infty,x_\infty)$ and we have
\begin{equation}\label{eq:stationary_compact}
    -\partial_a\rho_\infty + \lambda\partial_m(m\rho_\infty) - f(\varepsilon x_\infty)\rho_\infty + \delta_0^a(\gamma\circ\Pi)_*\left(f(\varepsilon x_\infty)\rho_\infty\right) = 0.
\end{equation}

Now, we write Eq.~\eqref{eq:PDE_compact_auxiliary} as the sum of a linear equation and a perturbation:
\begin{subequations} \label{eq:PDE_perturtation}
\begin{align}
    \partial_t(\rho_t, \zeta_t) &= \Lambda (\rho_t,\zeta_t) + (Z^{(1)}_t,Z^{(2)}_t), \\
    (\rho_0,\zeta_0) &= \left(u_0,0 \right),
\end{align}
where
\begin{align*}
    \Lambda(\rho_t,\zeta_t) &:= \left(-\partial_a\rho_t + \lambda\partial_m(m\rho_t) - f(\varepsilon x_\infty)\rho_t + \delta_0^a(\gamma\circ\Pi)_*\left(f(\varepsilon x_\infty)\rho_t\right), -\partial_s \zeta_t + \delta_0^s f(\varepsilon x_\infty)\rho_t\right),\\
    Z^{(1)}_t &:= [f(\varepsilon x_\infty) - f(\varepsilon x_t)]\rho_t + \delta_0^a(\gamma\circ\Pi)_*([f(\varepsilon x_t)-f(\varepsilon x_\infty)]\rho_t),\\
    Z^{(2)}_t &:= \delta_0^s[f(\varepsilon x_t)-f(\varepsilon x_\infty)]\rho_t.
\end{align*}
\end{subequations}
Let us put $\zeta_\infty(s) := f(\varepsilon x_\infty)\rho_\infty$, $\forall s\in\R_+$. Then, using Eq.~\eqref{eq:stationary_compact}, by the linearity of the operator $\Lambda$ and writing $\bar{\rho}_t := \rho_t-\rho_\infty$ and $\bar{\zeta}_t := \zeta_t-\zeta_\infty$, we get
\begin{subequations}\label{eq:PDE_perturbation_around}
\begin{align}
    \partial_t(\bar{\rho}_t, \bar{\zeta}_t) &= \Lambda (\bar{\rho}_t,\bar{\zeta}_t) + (Z^{(1)}_t,Z^{(2)}_t), \\
    (\bar{\rho}_0,\bar{\zeta}_0) &= \left(u_0-\rho_\infty,-\zeta_\infty \right).
\end{align}
\end{subequations}

Writing $(S^\Lambda_t)_{t\in\R_+}$ the semigroup associated with the operator $\Lambda$, we have, by Duhamel's formula,
\begin{equation}\label{eq:duhamel}
    (\bar{\rho}_t,\bar{\zeta}_t) = S^\Lambda_t(\bar{\rho}_0,\bar{\zeta}_0) + \int_0^t S^\Lambda_{t-s}(Z^{(1)}_s,Z^{(2)}_s)ds, \qquad \forall t\geq0.
\end{equation}

Let us define the weighted space
\begin{equation*}
    L^1_+(\mu) := \left\{\zeta \in L^1(\R_+\times\R_+\times\R_+^*,\R_+)\:\Big|\:\int_0^\infty \norm{\zeta(s)}_{L^1}\norm{h}_\infty e^{-\h s}ds<\infty\right\}.
\end{equation*}

Note that, for all $t\geq 0$,
\begin{align*}
    |x_t-x_\infty| &= \left|\int_0^\infty \int h(s)\zeta_t(s)\,dadm\,ds - \int_0^\infty \int h(s)\zeta_\infty(s)\,dadm\,ds\right| \\
    &\leq \int_0^\infty \norm{h}_\infty e^{-\h s}\norm{\zeta_t(s) - \zeta_\infty(s)}_{L^1}ds = \norm{\bar{\zeta}_t}_{L^1(\mu)}.
\end{align*}

Also, we have, for all $t\geq 0$,
\begin{subequations} \label{eq:bounds_Z}
\begin{align}
    \norm{Z^{(1)}_t}_{L^1} &\leq |\varepsilon| 2 L_f\norm{\rho_t}_{L^1}|x_t-x_\infty| \leq |\varepsilon| 2 L_f\norm{\bar{\zeta}_t}_{L^1(\mu)}, \\
    \norm{Z^{(1)}_t}_{L^1(w)} &\leq |\varepsilon| 2 L_f\norm{\rho_t}_{L^1(w)}|x_t-x_\infty| \leq |\varepsilon| 2 L_f\left(\norm{u_0}_{L^1(w)+\frac{b}{\alpha}}\right)\norm{\bar{\zeta}_t}_{L^1(\mu)}, \\
    \norm{Z^{(2)}_t}_{L^1(\mu)} &\leq |\varepsilon|\norm{h}_\infty L_f \norm{\rho_t}_{L^1}|x_t-x_\infty| \leq |\varepsilon|\norm{h}_\infty L_f \norm{\bar{\zeta}_t}_{L^1(\mu)},
\end{align}
\end{subequations}
where we have used Theorem~\ref{theorem:wellposedness}~$(ii)$ in the first line and Lemma~\ref{lemma:L1_stability} in the second.

\begin{lemma}\label{lemma:exp_decay_Lambda}
Grant Assumptions~\ref{assumption:functions} -- \ref{assumption:long-time} and \ref{assumption:nonlinear} and take $(\bar{\rho}_0,\bar{\zeta}_0)$ as in Eq.~\eqref{eq:PDE_perturbation_around}. 
There exists $K_1\geq1$ and $\mathfrak{a}_1>0$ such that, for all initial data $u_0\in L^1_+(w)$ with $\norm{u_0}_{L^1}=1$,
\begin{equation}\label{eq:exp_decay_Lambda}
    \norm{S^\Lambda_t(\bar{\rho}_0,\bar{\zeta}_0)}_{L^1(w)\times L^1(\mu)}\leq K_1 e^{-\mathfrak{a}_1 t}\norm{(\bar{\rho}_0, \bar{\zeta}_0)}_{L^1(w)\times L^1(\mu)}, \qquad \forall t\geq 0.
\end{equation}
If in addition, we grant Assumption~\ref{assumption:compact}, then there exists $K_2\geq1$ and $\mathfrak{a}_2>0$ such that, for all initial data $u_0\in L^1_+$ with $\norm{u_0}_{L^1}=1$,
\begin{equation}\label{eq:exp_decay_Lambda_doeblin}
    \norm{S^\Lambda_t(\bar{\rho}_0,\bar{\zeta}_0)}_{L^1\times L^1(\mu)}\leq K_2 e^{-\mathfrak{a}_2 t}\norm{(\bar{\rho}_0, \bar{\zeta}_0)}_{L^1\times L^1(\mu)}, \qquad \forall t\geq 0.
\end{equation}
\end{lemma}
\begin{proof}
We write $(S^\Lambda_t(\bar{\rho}_0,\bar{\zeta}_0)^{(1)}, S^\Lambda_t(\bar{\rho}_0,\bar{\zeta}_0)^{(2)}):= S^\Lambda_t(\bar{\rho}_0,\bar{\zeta}_0)$ the first and second component of $S^\Lambda_t(\bar{\rho}_0,\bar{\zeta}_0)$.

By Theorem~\ref{theorem:harris}, there exists $K\geq0$ and $\mathfrak{a}>0$ such that,
\begin{equation*}
    \norm{S^\Lambda_t(\bar{\rho}_0,\bar{\zeta}_0)^{(1)}}_{L^1(w)} \leq K e^{-\mathfrak{a} t}\norm{\bar{\rho}_0}_{L^1(w)}, \qquad \forall t\geq 0.
\end{equation*}
Then,
\begin{align*}
    \norm{S^\Lambda_t(\bar{\rho}_0,\bar{\zeta}_0)^{(2)}}_{L^1(\mu)} &= \int_0^t\norm{f(\varepsilon x_\infty)S^\Lambda_{t-s}(\bar{\rho}_0, \bar{\zeta}_0)^{(1)}}_{L^1}{C_h} e^{-\h s}ds \\
    &\qquad\qquad\qquad\qquad\qquad\qquad\qquad+ \int_t^\infty \norm{\bar{\zeta}_0(s)}_{L^1}{C_h} e^{-\h s}ds \\
    &\leq  {C_h} \left\{\norm{f}_\infty K \int_0^t e^{-\mathfrak{a}(t-s)}e^{-\h s}ds \norm{\bar{\rho}_0}_{L^1(w)} + e^{-\h t}\norm{\bar{\zeta}_0}_{L^1(\mu)}\right\}.
\end{align*}
Gathering the bounds on the two components and observing that the function $t\mapsto \int_0^t e^{-\mathfrak{a}(t-s)}e^{-\h s}ds$ decays exponentially, we conclude that there exists $K_1\geq 1$ and $\mathfrak{a}_1>0$ such that Eq.~\eqref{eq:exp_decay_Lambda} holds.

For Eq.~\eqref{eq:exp_decay_Lambda_doeblin}, we use Theorem~\ref{theorem:doeblin} and follow the same argument.
\end{proof}

We can now prove our main result:

\begin{proof}[Proof of Theorem~\ref{theorem:stability_nonlinear}]
By Duhamel's formula~\eqref{eq:duhamel}, Eq.~\eqref{eq:exp_decay_Lambda} in Lemma~\ref{lemma:exp_decay_Lambda} and the bounds Eqs.~\eqref{eq:bounds_Z}, for all $t\geq 0$,
\begin{align*}
    \norm{(\bar{\rho}_t,\bar{\zeta}_t)}_{L^1(w)\times L^1(\mu)} &\leq \norm{S^\Lambda_t(\bar{\rho}_0,\bar{\zeta}_0)}_{L^1(w)\times L^1(\mu)} + \int_0^t \norm{S^\Lambda_{t-s}(Z^{(1)}_s,Z^{(2)}_s)}_{L^1(w)\times L^1(\mu)}ds \\
    &\leq K_1 e^{-\mathfrak{a}_1 t}\norm{(\bar{\rho}_0,\bar{\zeta}_0)}_{L^1(w)\times L^1(\mu)} + K_1\int_0^t e^{-\mathfrak{a}_1 (t-s)}\norm{(Z^{(1)}_s,Z^{(2)}_s)}_{L^1(w)\times L^1(\mu)}ds \\
    &\leq K_1 e^{-\mathfrak{a}_1 t}\norm{(\bar{\rho}_0,\bar{\zeta}_0)}_{L^1(w)\times L^1(\mu)} +|\varepsilon|\tilde{C}_W\int_0^t e^{-\mathfrak{a}_1 (t-s)}\norm{(\bar{\rho}_s,\bar{\zeta}_s)}_{L^1(w)\times L^1(\mu)}ds\\
    &=: \mathcal{Q}(t),
\end{align*}
where $\tilde{C}_K$ is a constant depending on $W$. We have, for all $t\geq 0$,
\begin{align*}
    \frac{d}{dt}\mathcal{Q}(t) &= - \mathfrak{a}_1\mathcal{Q}(t) + |\varepsilon|\tilde{C}_W\norm{(\bar{\rho}_t,\bar{\zeta}_t)}_{L^1(w)\times L^1(\mu)} \\
    &\leq \left(-\mathfrak{a}_1 + |\varepsilon|\tilde{C}_W\right)\mathcal{Q}(t).
\end{align*}
Whence, by Grönwall's lemma,
\begin{equation*}
    \forall t\geq 0, \qquad \mathcal{Q}(t)\leq K_1 \norm{(\bar{\rho}_0,\bar{\zeta}_0)}_{L^1(w)\times L^1(\mu)}\exp\left(\left(-\mathfrak{a}_1 + |\varepsilon|\tilde{C}_W\right) t \right).
\end{equation*}
For all $t\geq 0$, we have
\begin{equation*}
    \norm{\rho_t - \rho_\infty}_{L^1(w)} + |x_t-x_\infty|\leq \norm{(\bar{\rho}_t,\bar{\zeta}_t)}_{L^1(w)\times L^1(\mu)} \leq \mathcal{Q}(t)
\end{equation*}
and 
\begin{equation*}
    \norm{\bar{\zeta}_0}_{L^1(\mu)}\leq \int_0^\infty \norm{f(\varepsilon x_\infty)\rho_\infty}_{L^1}{C_h} e^{-\h s}ds \leq \frac{\norm{f}_\infty {C_h}}{\h}.
\end{equation*}
Whence, choosing $\varepsilon^{**}_W := \frac{\mathfrak{a}_1}{\tilde{C}_W}\wedge \varepsilon^{*}$, we easily see that there exists $C\geq 1$ and $c_W>0$ such that Eq.~\eqref{eq:final} holds.

For Eq.~\eqref{eq:final_doeblin}, we use Eq.~\eqref{eq:exp_decay_Lambda_doeblin} instead of Eq.~\eqref{eq:exp_decay_Lambda} and follow the same argument.
\end{proof}

\section*{Appendix}
Here, we compare simulations of Eq.~\eqref{eq:PDE_aSRM0} with simulations of the time elapsed neuron network model \cite{PakPer09}.

If, the firing rate function $f$ does not depend on $m$ and if we put
\begin{equation}\label{eq:f_TENNM}
f(a, \varepsilon x_t):=\hat{f}(\eta(a) + \varepsilon x_t), 
\end{equation}
then Eq.~\eqref{eq:PDE_aSRM0} reduces to the time elapsed neuron network model
\begin{subequations} \label{eq:TENNM}
\begin{align} 
    &\partial_t\rho_t(a)+\partial_a\rho_t(a)= - f(a, \varepsilon x_t)\rho_t(a), \\
    &\rho_t(0) =\int_0^\infty f(a,\varepsilon x_t)\rho_t(a)da,  \\
    &x_t = \int_0^t  h(t-s)\int_{0}^\infty\int_0^\infty f(a,\varepsilon x_s)\rho_s(a)dads,\\
    &\rho_0(a) = u_0(a).
\end{align}
\end{subequations}
Eq.~\eqref{eq:TENNM} is the population equation for non-adaptive $\text{SRM}_0$ neurons (or age-dependent nonlinear Hawkes processes) \cite{Che17}. As reported previously, Eq.~\eqref{eq:TENNM} exhibits self-sustained oscillations for large $\varepsilon$ or relaxation to a stationary state for small $\varepsilon$ (see Fig.~2).
\begin{figure} \label{fig:2}
  \begin{center}
\includegraphics[scale=0.65]{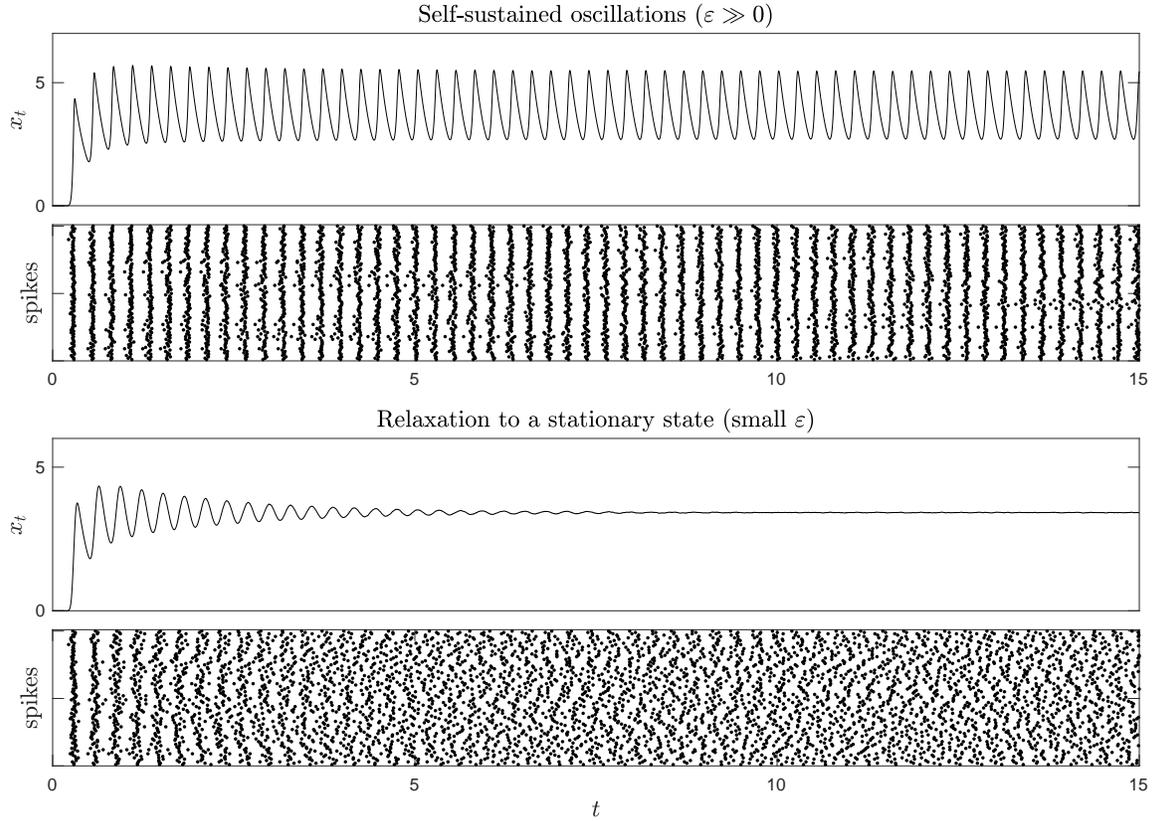}
\caption{\textbf{Same as Fig.~1 but for the time elasped neuron network model Eq.~\eqref{eq:TENNM}}. Simulations of a network of $5\cdot 10^5$ non-adaptive $\text{SRM}_0$ neurons, approximating Eq.~\eqref{eq:TENNM}, with identical parameters (except for $\varepsilon$) and identical initial conditions. Neuronal parameters are the same is in Fig.~1, expect that $f$ is replaced by Eq.~\eqref{eq:f_TENNM}. The $\varepsilon$ have also been adapted.}
  \end{center}
\end{figure}

\section*{Acknowledgements}
We thank Stéphane Mischler for supervising this work and Wulfram Gerstner for his comments on the manuscript. This research has been funded by the Swiss National Science Foundation (grant no.~200020\_184615) and the European Union's Horizon 2020 research and innovation programme under the Marie Skłodowska-Curie grant agreement No~754362. \includegraphics[scale=0.15]{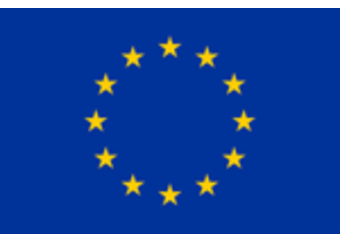}

\bibliographystyle{plain} 
\bibliography{mybib}
\end{document}